\documentclass[10pt]{article}
\usepackage[dvips]{graphicx}
\usepackage{subfigure}
\usepackage{pdfpages}
\usepackage[english]{babel}
\usepackage{blindtext}
\usepackage{algorithm} 
\usepackage{algpseudocode}
\usepackage{multicol}
\usepackage{amsmath}
\usepackage[morefloats=18]{morefloats}
\usepackage{bm}
\usepackage[figuresright]{rotating}
\usepackage{threeparttable}
\usepackage{dsfont}
\usepackage{amsfonts}
\usepackage{booktabs}
\usepackage{cite}
\usepackage{rotating}
\usepackage{amssymb}
\usepackage{float}
\usepackage{color}
\usepackage{algorithm}
\usepackage{algpseudocode}
\usepackage{amsmath}
\usepackage{amsthm}

\newtheorem{theorem}{Theorem}[section]
\newtheorem{remark}{Remark}[section]

\usepackage{CJK}
\begin{document}

\begin{center}
{\Large \bf Generalized Golub-Kahan bidiagonalization for generalized saddle point systems 
}\\
[0.2cm]{\normalsize Na-Na Wang$^{a}$\footnote{Corresponding author.\\
\hspace*{0.5cm} Email:~wangnana@chd.edu.cn (Na-Na Wang).}
\quad\quad Ji-Cheng Li$^{b}$
}
\\
[0.2cm] {\normalsize $a$. School of Science,
        Chang'an University, Xi'an, Shaanxi, 710064, China}
\\
[0.2cm] {\normalsize $b$. School of Mathematics and Statistics,
        Xi'an Jiaotong University, Xi'an, Shaanxi, 710049, China}
\end{center}

\vspace{0.5cm}

\large

\date{}



{\bf\large Abstract}: We consider the iterative solution of generalized saddle point systems. When the right bottom block is zero, Arioli [SIAM J. Matrix Anal. Appl., 34 (2013), pp. 571--592] proposed a CRAIG algorithm based on generalized Golub-Kahan Bidiagonalization (GKB) for the augmented systems with the leading block being symmetric and positive definite (SPD), and then Dumitrasc et al. [SIAM J. Matrix Anal. Appl., 46 (2025), pp. 370--392] extended the GKB for the case where the symmetry condition of the leading block no longer holds and then proposed nonsymmetric version of the CRAIG (nsCRAIG) algorithm. The CRAIG and nsCRAIG algorithms are theoretically equivalent to the Schur complement reduction (SCR) methods where the Conjugate Gradient (CG) method and the Full Orthogonalization Method (FOM) are applied to the associated Schur-complement equation, respectively.
We extend the GKB and its nonsymmetric counterpart 
used separately in CRAIG and nsCRAIG algorithms for the case where the right bottom block of saddle point system is nonzero. On this basis,
we propose CRAIG and nsCRAIG algorithms for the solution of the generalized saddle point problems with the leading block being SPD and nonsymmetric positive definite (NSPD), respectively. They are also theoretically equivalent to the SCR methods with inner CG and FOM iterations for the associated Schur-complement equation, respectively.
Moreover, we give algorithm steps of the two new solvers and propose appropriate stopping criteria based on an estimate of the energy norm for the error and the residual norm. 
Numerical comparison with MINRES or GMRES highlights the advantages of our proposed strategies regarding its high computational efficiency and/or low memory requirements and the associated implications.

{\bf AMS(2000)}: \quad 65F10; \quad 65F50; \quad ‌65N22


\vspace*{0.5cm}

\numberwithin{equation}{section}
\section{Introduction}

\quad We consider the iterative solution of the following generalized saddle point problem
\begin{equation}\label{e11}
\mathcal{A}z \triangleq \left[\begin{array}{c c}
		M & A\\
		A^T & -C\\
	\end{array}\right]
	\left[\begin{array}{c}
		u \\ p
	\end{array}\right]
	=\left[\begin{array}{c}
		0 \\ b
	\end{array}\right] \triangleq f,
\end{equation}
where $M\in \mathbb{R}^{m \times m}$ is SPD or NSPD,
$A\in \mathbb{R}^{m \times n}$ has full column rank, and $C\in \mathbb{R}^{n \times n}$ is nonzero and symmetric positive semi-definite (SPSD), $b\in \mathbb{R}^{n}$,
and $n\leq m$. 
For a generalized saddle point problem with a general right-hand side vector:
\begin{equation}\label{e12}
\left[\begin{array}{c c}
		M & A\\
		A^T & -C\\
	\end{array}\right]
	\left[\begin{array}{c}
		w \\ p
	\end{array}\right]
	=\left[\begin{array}{c}
		b_1 \\ b_2
	\end{array}\right],
\end{equation}
with $b_1\in \mathbb{R}^{m}$ and $b_2\in \mathbb{R}^{n}$, by using the transformation of the form
\begin{equation*}
w_0=M^{-1}b_1, \quad u=w-w_0, \quad b=b_2-A^Tw_0,
\end{equation*}
the resulting equivalent system is (\ref{e11}). Once $u$ has been computed, $w$ can be recovered as $w=u+w_0$.
In the following, we always use ``saddle point problem" and ``generalized saddle point problem" to denote the linear system (\ref{e11}) with $C=O$ and $C\neq O$, respectively, and use ``symmetric" and ``nonsymmetric" (generalized) saddle point problem to denote the linear system (\ref{e11}) with $A$ SPD and NSPD, respectively.

The linear systems (\ref{e11}) with nonsymmetric $M$ arise in the stabilized finite-element discretization of Oseen problems obtained by linearization, through Picard's method, of the steady-state Navier-Stokes equations governing the flow of a Newtonian incompressible viscous fluid \cite{MBGHGJL2005}.
The system (\ref{e11}) with symmetric $M$ comes from the stabilized finite-element discretization of the steady-state Stokes equations governing the flow of a slow-moving, highly viscous fluid.
When both $M$ and $C$ are SPD, the system (\ref{e11}) is called symmetric quasi-definite system, which may be interpreted as regularized linear least-squares problem in appropriate metrics. The type of system originates from applications such as regularized interior-point methods for convex optimization and stabilized control problems \cite{MBGHGJL2005,MDAVDSDDS2010,MPFDO2012,DOMA2017} and preconditioned iterative methods for ill-posed constrained and weighted toeplitz least squares problems with Tikhonov regularization \cite{MBMKN2006}.
The systems (\ref{e11}) arise also from Lagrangian approaches for variational problems with equality constraints when the constraints are relaxed or a penalty term is applied \cite{JPAJW2015}. In the case, $M$ is usually symmetric but may also be nonsymmetric and often has additional properties, e.g., it accounts for local convexity of the optimization problem.

The solution algorithms for generalized saddle point problems (\ref{e11}) can be subdivided into two broad categories, called coupled (or `all-at-once') and segregated methods \cite{MBGHGJL2005}.
Coupled methods deal with the system (\ref{e11}) as a whole, computing $u$ and $p$ (or approximations to them) simultaneously. 
These methods include both direct solvers based on triangular factorizations of the global coefficient matrix of (\ref{e11}), and iterative algorithms like Krylov subspace methods \cite{YS2003} applied to the entire system (\ref{e11}), typically with some form of preconditioning.

Segregated methods, on the other hand, involve the solution of two linear systems of size smaller than $m + n$ (called reduced systems), one for each of $u$ and $p$. Segregated methods can be either direct or iterative, or involve a combination of the two. The main representatives of the segregated approach are the null space method \cite{HSDNIMGWHASAJW2006,DDSDO2021}, which relies on a basis for the null space for the constraints and the SCR 
method \cite{MBGHGJL2005}, which is based on a block LU factorization of the generalized saddle point matrix in (\ref{e11}).
Specifically, 
the SCR method for  (\ref{e11}) is introduced as follows. By the block LU factorization of the coefficient matrix of (\ref{e11}):
\begin{equation*}
\left[\begin{array}{c c}
		M & A\\
		A^T & -C\\
	\end{array}\right]=
\left[\begin{array}{c c}
		I_m & O \\
		A^TM^{-1} & I_n\\
	\end{array}\right]\left[\begin{array}{c c}
		M & A\\
		O & -S\\
	\end{array}\right],
\end{equation*}
with $S=A^TM^{-1}A+C$ being the Schur complement of the (1,1)-block $M$ in (\ref{e11}), we get the transformed system
\begin{equation*}
\left[\begin{array}{c c}
		M & A\\
		O & -S\\
	\end{array}\right]
	\left[\begin{array}{c}
		u \\ p
	\end{array}\right]
	=\left[\begin{array}{c c}
		I_m & O \\
		-A^TM^{-1} & I_n\\
	\end{array}\right]\left[\begin{array}{c}
		0 \\ b
	\end{array}\right]=\left[\begin{array}{c}
		0 \\ b
	\end{array}\right].
\end{equation*}
Solving this block upper triangular system by block backsubstitution leads to the two reduced systems of order $n$ for $p$ and $m$ for $u$:
\begin{equation}\label{e13}
(A^TM^{-1}A+C)p=Sp=-b,
\end{equation}
\begin{equation}\label{e14}
Mu=-Ap.
\end{equation}
Once $p_{*}$ has been computed from (\ref{e13}), $u_{*}$ can be obtained by solving (\ref{e14}). These systems can be solved either directly or iteratively.
In this paper, we interested in the development of some segregated methods based on the Golub-Kahan bidiagonalization and its nonsymmetric variant for generalized saddle point problems (\ref{e11}).

The generalized CRAIG solver or the GKB algorithm for symmetric saddle point problems (\ref{e11}) 
was introduced by Arioli in \cite{MA2013} and it is based on the bidiagonalization of the (1,2)-block $A$ of the system matrix.
This solver 
is theoretically equivalent to the SCR method where the CG iteration is applied to the associated SPD Schur-complement system (\ref{e13}) with an SPD preconditioner and then the direct method such as Cholesky factorization is applied to the SPD system (\ref{e14}). The latter is abbreviated as SCR(CG).
The GKB was further extended to the (1,2)-block $A$ of the \textbf{nonsymmetric} saddle point problem (\ref{e11}), 
and on this basis, 
a nonsymmetric version of the CRAIG algorithm, called nsCRAIG, was introduced by Dumitrasc et al. in \cite{ADCKUR2025}.
The nsCRAIG algorithm 
is theoretically equivalent to the SCR method where the FOM method is applied to the related NSPD Schur-complement system (\ref{e13}) and then the direct method such as LU factorization is applied to the NSPD system (\ref{e14}). The latter is abbreviated as SCR(FOM).
It is worth to emphasize that to apply CRAIG and nsCRAIG to the saddle point problems (\ref{e12}), the information in the right-hand side needs to compress into the lower block of the right-hand side as in (\ref{e11}), acts as a starting point for the bidiagonalization in CRAIG and the decomposition in nsCRAIG, see the next section. There are works \cite{DOMA2017,SMA1995} that use another version of the GKB, where the compression leads to $[b;0]$ with $b \in R^m$.

The CRAIG for symmetric saddle point problem \cite{MA2013} and LSQR and LMSR for Least-Squares \cite{CCPMAS1982,DCLFMS2011}, are examples of solvers stemming from Golub-Kahan bidiagonalization. Their use can be preferable to the alternative of using CG directly on the Schur-complement equation or normal equation, which would involve operating with a squared condition number and the difficulty in ensuring the accuracy of its solution, as described in
\cite[Chapter 8.1]{YS2003}, and \cite{CCPMAS1982}. With the similar reason, the use of the recently proposed nsCRAIG for nonsymmetric saddle point problem \cite{ADCKUR2025}, an example of solvers stemming from the nonsymmetric version of GKB, can be preferable to using FOM directly on the Schur-complement equation. Therefore, the above advantages of the CRAIG and nsCRAIG for saddle point problem motivate us to extend these algorithms to solve the generalized saddle point problems (1.1) and obtain two new solvers. Compared to CRAIG and nsCRAIG for the saddle point problems \cite{MA2013,ADCKUR2025}, the obtained two new algorithms are based on the bidiagonalization and decomposition of the \textbf{augmentation of} (1,2)-block of (\ref{e11}), respectively, rather than the bidiagonalization and decomposition of the (1,2)-block of (\ref{e11}).
Theoretically, the two new solvers are also equivalent to the SCR(CG) and SCR(FOM) methods, respectively.

By exploiting a suitable reformulation of (\ref{e11}) suggested by Dollar et al. \cite{HSDNIMGWHASAJW2006}, see also \cite{DOMA2017,DDSDO2021}, we reformulate it 
as a system with saddle point matrix of zero (2,2)-block as follows.
Assume that $\mbox{rank}(C) = l$ and $C$ has been decomposed as
\begin{equation}\label{e15}
C=E^TFE,
\end{equation}
where $F\in \mathbb{R}^{l \times l}$ is SPD since $C$ is SPSD and $E\in \mathbb{R}^{l \times n}$. Then, by using the auxiliary variable
\begin{equation*}
c=-FEp,
\end{equation*}
the system (\ref{e11}) may be rewritten as
\begin{equation}\label{e16}
\left[\begin{array}{ccc}
		M && A \\
		&F^{-1}& E\\
		A^T &E^T & \\
	\end{array}\right]
	\left[\begin{array}{c}
		u \\ c \\ p
	\end{array}\right]
	=\left[\begin{array}{c}
		0 \\ 0 \\ b
	\end{array}\right],
\end{equation}
which has a standard saddle point form
\begin{equation}\label{e17}
\left[\begin{array}{c c}
		\bar{M} & \bar{A}\\
		\bar{A}^T & \\
	\end{array}\right]
	\left[\begin{array}{c}
		\bar{u} \\ p
	\end{array}\right]
	=\left[\begin{array}{c}
		0 \\ b
	\end{array}\right],\quad
\bar{M} = \left[\begin{array}{c c}
		M & \\
		   & F^{-1} \\
	\end{array}\right], \quad \bar{A}=\left[\begin{array}{c}
		A \\ E
	\end{array}\right], \quad
\bar{u} = \left[\begin{array}{c}
		u \\ c \\
	\end{array}\right],
\end{equation}
where $\bar{M} \in R^{(m+l)\times (m+l)}$ is SPD if $M$ is SPD and NSPD if $M$ is NSPD,
and $\bar{A}\in R^{(m+l)\times n}$ has full column rank since 
$A$ has full column rank. Besides (\ref{e16}), the system (\ref{e11}) has another equivalent saddle point form
\begin{equation}\label{e18}
\left[\begin{array}{ccc}
		M && A\\
		&I_n& C\\
		A^T &I_n & \\
	\end{array}\right]
	\left[\begin{array}{c}
		u \\ a \\ p
	\end{array}\right]
	=\left[\begin{array}{c}
		0 \\ 0 \\ b
	\end{array}\right],
\end{equation}
where $a=-Cp$. If $C\neq I_n$, then the system (\ref{e18}) is nonsymmetric whether $M$ is symmetric or not. However, its Schur complement $A^TM^{-1}A+C$ is SPD if $M$ is SPD. In addition, the size of (\ref{e18}) is larger than that of (\ref{e16}) if $l<n$.
Obviously, once the upper and lower blocks $u$ and $p$ of the solution for system (\ref{e18}) have been obtained, then the solution for system (\ref{e11}) is obtained.


The organization of this paper is as follows. In subsections \ref{SubSec21} and \ref{SubSec22} of section \ref{Sec2}, we separately take the symmetric and nonsymmetric saddle point problems (\ref{e17}) as examples to introduce the existing bidiagonalization and decomposition of the (1,2)-block of the systems.
In section \ref{Sec3}, we extend the considerations in subsection \ref{SubSec21} to the symmetric generalized saddle point system (\ref{e11}) and develop the bidiagonalization of the augmentation of the (1,2)-block of the system, and then use it to derive the CRAIG solver and its stopping criteria.
In section \ref{Sec4}, we extend the considerations in subsection \ref{SubSec22} to the nonsymmetric generalized saddle point system (\ref{e11}) and introduce an adapted decomposition of the augmentation of the (1,2)-block of the system, the corresponding nsCRAIG solver and its stopping criteria.
In section \ref{Sec5}, 
numerical experiments on some Stokes equations and Navier-Stokes equations are used to verify that the numerically equivalence between the proposed CRAIG and nsCRAIG solvers and the SCR method with inner CG and FOM iterations, respectively,
and the advantages of our proposed solvers over the common methods for generalized saddle point problems, i.e., the MINRES and GMRES with an appropriate block diagonal preconditioner.
In section \ref{Sec6}, some conclusions are given.

For simplicity of description, some notations and assumptions are presented here. For two symmetric matrices $H_1$ and $H_2$, $H_1\succ H_2$ ($H_1\succeq H_2$) means that $H_1-H_2$ is SPD (SPSD).
For any matrix $H$, $\mbox{Null}(H)$ and $\mbox{Range}(H)$ represent its null space and range space, respectively. The zero matrix $O$ and the identity matrix $I$ or $I_{k}$ are used according to the appropriate dimensions with $k$ being any positive integer.
For vectors $x\in R^m$ and $y\in R^n$, $[x; y] \in R^{m+n}$ denotes a column vector like in the Matlab.

\section{The GKB for symmetric and nonsymmetric saddle point problems}\label{Sec2}

\quad In this section, we review the existing GKB for the symmetric saddle point problems of the form (\ref{e17}), and its connections with the tridiagonalization of the related Schur complement. 
Then, we review the existing nonsymmetric version of GKB for the nonsymmetric systems (\ref{e17}) and its connections with the upper Hessenberg form of the associated Schur complement.


\subsection{The GKB for symmetric saddle point problems (\ref{e17})}\label{SubSec21}

\quad In this subsection, we review the GKB as building block of the generalized CRAIG algorithm introduced in \cite{MA2013,VDADCKUR2023,ADCKUR2025} for the solution of the symmetric saddle point systems (\ref{e17}). 


Let $N\in R^{m\times m}$ be an SPD matrix. To properly describe the GKB, we need to define the following inner product and norms
\begin{equation}\label{e21}
\begin{array}{c}
\langle u,v \rangle_{M}=u^TMv,\quad \|u\|_{M}=\sqrt{u^TMu},\quad \langle c,d \rangle_{F^{-1}}=c^TF^{-1}d,\quad \|c\|_{F^{-1}}=\sqrt{c^TF^{-1}c},\\
\langle \bar{u},\bar{v} \rangle_{\bar{M}} = \bar{u}^T\bar{M}\bar{v}= 
\langle u,v \rangle_{M}+\langle c,d \rangle_{F^{-1}},\quad 
\|\bar{u}\|_{\bar{M}}=\sqrt{\bar{u}^T\bar{M}\bar{u}}
=\sqrt{\|u\|_M^2+\|c\|_{F^{-1}}^2},\\
\langle x,y \rangle_{N}=x^TNy,\quad \|y\|_{N}=\sqrt{y^TNy},  \quad  \|y\|_{N^{-1}}=\sqrt{y^TN^{-1}y},
\end{array}
\end{equation}
where $\bar{u}=[u;c]$ and $\bar{v}=[v;d]$ with $u,v\in R^m$ and  $c,d\in R^l$, and $x,y\in R^n$.

Given the right-hand side vector $b \in R^{n}$, the first step of the generalized GKB 
of $\bar{A}$ is 
\begin{equation}\label{e22}
\beta_1=\|b\|_{N^{-1}}, \quad q_1=N^{-1}b/\beta_1.
\end{equation}
After $k$ iterations, the partial bidiagonalization of $\bar{A}$ is iteratively given by
\begin{equation}\label{e23}
\left\{
\begin{array}{ll}
\bar{A}Q_k=\bar{M}\bar{V}_kB_k, & \bar{V}_k^T\bar{M}\bar{V}_k=I_k,\\
\bar{A}^T\bar{V}_k=NQ_kB_k^T+\beta_{k+1}Nq_{k+1}e_k^T, & Q_k^TNQ_k=I_k,
\end{array}
\right.
\end{equation}
with the upper bidiagonal matrix
\begin{equation}\label{e24}
B_k=\left[
\begin{array}{ccccc}
\alpha_1 & \beta_2 & 0 & \dots & 0\\
0& \alpha_2 & \beta_3 & \dots & 0\\
\vdots & \vdots  & \vdots  & \vdots  & \vdots  \\
0 & \dots & 0& \alpha_{k-1} & \beta_k \\
0 & \dots & 0& 0 & \alpha_{k} \\
\end{array}
\right]\in R^{k\times k},
\end{equation}
$\bar{V}_k=[\bar{v}_{1},\bar{v}_{2},\dots, \bar{v}_{k}] \in R^{(m+l)\times k}$ and $Q_k=[q_{1},q_{2},\dots, q_{k}] \in R^{n\times k}$. The basis $\bar{V}_k$ has $\bar{M}$-orthogonal columns with respect to (w.r.t.) the inner product $\langle \cdot, \cdot \rangle_{\bar{M}}$ and norm $\| \cdot\|_{\bar{M}}$.
Similarly, the basis $Q_k$ has $N$-orthogonal columns w.r.t. the inner product $\langle  \cdot,  \cdot \rangle_{N}$ and norm $\| \cdot\|_{N}$.
It is sufficient to store only the latest left vector $\bar{v}_{k}$ and use it to compute $\bar{v}_{k+1}$. The same is true for the right vector $q_{k+1}$.
Prior to the normalization leading to $\bar{v}_{k+1}$ and $q_{k+1}$, the norms are stored as $\alpha_{k+1}$ for $\bar{v}_{k+1}$ and $\beta_{k+1}$ for $q_{k+1}$, as detailed in Algorithm \ref{alg_S_GKB},
which can be used to define CRAIG solver \cite{MA2013,VDADCKUR2023}.
In the algorithm, $\bar{v}_k=[v_{k,x};v_{k,c}]$ and $\bar{w}_k=[w_{k,x};w_{k,c}]$ with $v_{k,x}, w_{k,x} \in R^m$ and $v_{k,c}, w_{k,c} \in R^l$.

\begin{algorithm}[!ht]
\caption{Golub-Kahan bidiagonalization \cite{MA2013,VDADCKUR2023}}
\label{alg_S_GKB}
\textbf{Require}: $M\in R^{m\times m}$ SPD, $F\in R^{l\times l}$ SPD, $A\in R^{m\times n}$ with full column rank, $E\in R^{l\times n}$, $b\in R^{n}$, maxit
\\
\textbf{Output}: $q_{k}$, $\beta_k$, $v_{k,x}$, $v_{k,c}$, $\alpha_k$, ($k=1,2,\cdots,$ maxit)
  \begin{algorithmic}[1]
    \State $\beta_1=\|b\|_{N^{-1}}$; $q_1=N^{-1}b/\beta_1$
    \State $w_{1,x} = M^{-1}Aq_1$; $w_{1,c} = FEq_1$
    {\hfill \textcolor{gray}{$\bar{w}_1 = \bar{M}^{-1}\bar{A}q_1$}}
    \State $\alpha_1 = \sqrt{\|w_{1,x}\|_{M}^2+\|w_{1,c}\|_{F^{-1}}^2}$
    {\hfill \textcolor{gray}{$\alpha_1 = \|\bar{w}_1\|_{\bar{M}}$}}
    \State $v_{1,x} = w_{1,x}/\alpha_1$; $v_{1,c} = w_{1,c}/\alpha_1$
    {\hfill \textcolor{gray}{$\bar{v}_1 = \bar{w}_1/\alpha_1$}}
    \State $k=1$
    \While{$k<$ maxit}
      \State $g_k=N^{-1}(A^Tv_{k,x}+E^Tv_{k,c}-\alpha_kNq_k)$
      {\hfill \textcolor{gray}{$g_k=N^{-1}(\bar{A}^T\bar{v}_k-\alpha_kNq_k)$}}
      \State $\beta_{k+1}=\|g_k\|_{N}$; $q_{k+1}=g_k/\beta_{k+1}$
      \State $w_{k+1,x}=M^{-1}(Aq_{k+1}-\beta_{k+1}Mv_{k,x})$ 
      {\hfill \textcolor{gray}{$\bar{w}_{k+1}=\bar{M}^{-1}(\bar{A}q_{k+1}-\beta_{k+1}\bar{M}\bar{v}_k)$}}
      \State $w_{k+1,c}=F(Eq_{k+1}-\beta_{k+1}F^{-1}v_{k,c})$
      \State $\alpha_{k+1}=\sqrt{\|w_{k+1,x}\|_{M}^2+\|w_{k+1,c}\|_{F^{-1}}^2}$
      {\hfill \textcolor{gray}{$\alpha_{k+1}=\|\bar{w}_{k+1}\|_{\bar{M}}$}}
      \State $v_{k+1,x}=w_{k+1,x}/\alpha_{k+1}$; $v_{k+1,c}=w_{k+1,c}/\alpha_{k+1}$
      {\hfill \textcolor{gray}{$\bar{v}_{k+1}=\bar{w}_{k+1}/\alpha_{k+1}$}}
      \State $k=k+1$
    \EndWhile
  \end{algorithmic}
\end{algorithm}

From (\ref{e23}), we can draw a link between the bidiagonalization of the (1,2)-block $\bar{A}$ of (\ref{e17}) and the tridiagonalization of the centered-preconditioned Schur complement $N^{-\frac{1}{2}}SN^{-\frac{1}{2}}$ with $S=\bar{A}^T\bar{M}^{-1}\bar{A}$:
\begin{equation}\label{e25}
\begin{array}{ll}
N^{-\frac{1}{2}}SN^{-\frac{1}{2}}N^{\frac{1}{2}}Q_k & =N^{-\frac{1}{2}}\bar{A}^T\bar{M}^{-1}\bar{A}Q_k=N^{-\frac{1}{2}}\bar{A}^T\bar{V}_kB_k \\
& =N^{-\frac{1}{2}}(NQ_kB_k^T+\beta_{k+1}Nq_{k+1}e_k^T)B_k \\
& =N^{\frac{1}{2}}Q_k(B_k^TB_k)+\alpha_k\beta_{k+1}N^{\frac{1}{2}}q_{k+1}e_k^T, \\
\end{array}
\end{equation}
where $B_k^TB_k$ is the tridiagonal matrix given by the Lanczos process specific to CG. In the last line of (\ref{e25}), we use the fact that $e_k^TB_k=\alpha_ke_k^T$ by the structure of $B_k$ in (\ref{e24}).
When $k=n$, the bidiagonalization residual $\beta_{k+1}Nq_{k+1}e_k^T$ in (\ref{e23}) vanishes and we obtain the tridiagonalization of  $N^{-\frac{1}{2}}SN^{-\frac{1}{2}}=(N^{\frac{1}{2}}Q_n)(B_n^TB_n)(N^{\frac{1}{2}}Q_n)^T$ by (\ref{e25}).
Therefore, we can implicitly tridiagonalize $N^{-\frac{1}{2}}SN^{-\frac{1}{2}}$ in the sense that we find $Q_n$ and $B_n$ without knowledge of $S$.
Since $N^{-\frac{1}{2}}SN^{-\frac{1}{2}}$ is symmetric, the column vectors in basis $Q_k$ generated in Lanczos process 
satisfy the short recurrences, where one vector needs to be orthogonalized only against the previous two w.r.t. the inner product $\langle \cdot, \cdot \rangle_{N}$. Then, in exact arithmetic, it will also be orthogonal to all the previous ones w.r.t. $\langle  \cdot,  \cdot \rangle_{N}$. Consequently, only the most recent vectors are stored and used.

Let 
\begin{equation}\label{e26}
\bar{V}_{k}=[V_{k,x}^T, V_{k,c}^T]^T,
\quad V_{k,x}=[v_{1,x},v_{2,x},\dots, v_{k,x}]\in R^{m \times k}, \quad V_{k,c}=[v_{1,c},v_{2,c},\dots, v_{k,c}]\in R^{l \times k}.
\end{equation}
Then,
by (\ref{e17}) and (\ref{e26}), the (\ref{e23}) can be rewritten as
\begin{equation}\label{e27}
\left\{
\begin{array}{ll}
AQ_k=MV_{k,x}B_k,\quad EQ_k=F^{-1}V_{k,c}B_k, & V_{k,x}^TMV_{k,x}+V_{k,c}^TF^{-1}V_{k,c}=I_k,\\
A^TV_{k,x}+E^TV_{k,c}=NQ_kB_k^T+\beta_{k+1}Nq_{k+1}e_k^T, & Q_k^TNQ_k=I_k.
\end{array}
\right.
\end{equation}

In Section \ref{Sec3}, we show how the GKB of $\bar{A}$ in (\ref{e23}) or
(\ref{e27})
can be further reformulated by referring to the original generalized saddle point system (\ref{e11}), thus avoiding the use of $E$ and $F$ and the factorization (\ref{e15}).

\subsection{The GKB for nonsymmetric saddle point problems (\ref{e17})}\label{SubSec22}


\quad In \cite{ADCKUR2025}, the nonsymmetric version of GKB as building block of nsCRAIG algorithm has been introduced. In this subsection, corresponding to the symmetric case, we provide the preconditioned version of the nonsymmetric GKB for the solution of the nonsymmetric saddle point systems (\ref{e17}).



Given the right-hand side vector $b \in R^{n}$, the first step of the decomposition of $\bar{A}$ is identical to that of the GBK given in (\ref{e22}).
After $k$ iterations, the partial decomposition of $\bar{A}$ is iteratively given by \begin{equation}\label{e28}
\left\{
\begin{array}{ll}
\bar{A}Q_k=\bar{M}\bar{V}_kB_k, & \bar{V}_k^T\bar{M}\bar{V}_k=L_k,\\
\bar{A}^T\bar{V}_k=NQ_kH_k+\beta_{k+1}Nq_{k+1}e_k^T, & Q_k^TNQ_k=I_k,
\end{array}
\right.
\end{equation}
where the upper bidiagonal matrix $B_k\in R^{k\times k}$ is given as in (\ref{e24}) and
\begin{equation}\label{e29}
H_k =\left[
\begin{array}{ccccc}
h_{1,1} & h_{1,2} & h_{1,3} & \dots & h_{1,k}\\
\beta_2& h_{2,2} & h_{2,3} & \dots & h_{2,k}\\
0 & \beta_3  & h_{3,3}  & \dots  & h_{3,k} \\
\vdots & \vdots & \vdots & \vdots & \vdots \\
0 & 0 & 0& \beta_k & h_{k,k} \\
\end{array}
\right]\in R^{k\times k},
\end{equation}
$Q_k\in R^{n\times k}$ is defined as in (\ref{e23}), $\bar{V}_k\in R^{(m+l)\times k}$ and $L_k\in R^{k\times k}$.
By \cite[Proposition 3.1]{ADCKUR2025}, we know that the matrix $L_k$ is unit lower triangular and $H_k=B_k^TL_k^T$ that can be considered as the LU decomposition of $H_k$.
Similar to the symmetric case, the basis $Q_k$ has $N$-orthogonal columns w.r.t. the inner product $\langle  \cdot,  \cdot \rangle_{N}$ and norm $\| \cdot\|_{N}$. However, in the case that where $\bar{M}$ is not symmetric but positive definite since $M$ is NSPD, $\bar{M}$ itself cannot define a norm, but its symmetric part can define a norm well, which can still be written in the form defined in (\ref{e21}). At this point, the orthogonality of $\bar{M}$ between the columns of the basis $\bar{V}_k$ can be intuitively considered a kind of one-sided, nonsymmetric orthogonality, see \cite{ADCKUR2025}.

By using the notations similar to that in Algorithm \ref{alg_S_GKB}, we give the nonsymmetric version of GKB in Algorithm \ref{alg_NS_GKB}, which can be used to define nsCRAIG solver. When $N=I_n$, Algorithm \ref{alg_NS_GKB} reduces the nonsymmetric version of GKB as building block of nsCRAIG proposed in \cite{ADCKUR2025}.
In any step $k> 0$ of Algorithm \ref{alg_NS_GKB}, it is necessary to store all the right vectors in $Q_{k}$ and use them in the orthogonalization process to maintain global mutual orthogonality w.r.t. $\langle  \cdot,  \cdot \rangle_{N}$, where the vector $h_k = [h_{1,k}, h_{2,k}, \cdots, h_{k,k}]^T$ is computed to store the inner products $\langle  \cdot,  \cdot \rangle_{N}$ between $\hat{g}_k$ and all the columns in $Q_{k}$. It is sufficient to store only the latest left vector $\bar{v}_{k}$ and use it to compute $\bar{v}_{k+1}$. Prior to the normalization leading to $\bar{v}_{k+1}$ and $q_{k+1}$, the norms are stored as $\alpha_{k+1}$ for $\bar{v}_{k+1}$ and $\beta_{k+1}$ for $q_{k+1}$, as detailed in Algorithm \ref{alg_NS_GKB}, similar to Algorithm \ref{alg_S_GKB}.

\begin{algorithm}[!ht]
\caption{The nonsymmetric version of GKB \cite{ADCKUR2025}}
\label{alg_NS_GKB}
\textbf{Require}: $M\in R^{m\times m}$ NSPD, $F\in R^{l\times l}$ SPD, $A\in R^{m\times n}$ with full column rank, $E\in R^{l\times n}$, $b\in R^{n}$, maxit
\\
\textbf{Output}: $q_{k}$, $\beta_k$, $v_{k,x}$, $v_{k,c}$, $\alpha_k$, ($k=1,2,\cdots,$ maxit), $h_k$, ($k=1,2,\cdots,$ maxit$-1$)
  \begin{algorithmic}[1]
    \State $\beta_1=\|b\|_{N^{-1}}$; $q_1=N^{-1}b/\beta_1$; $Q_{1}=q_{1}$
    \State $w_{1,x} = M^{-1}Aq_1$; $w_{1,c} = FEq_1$
    {\hfill \textcolor{gray}{$\bar{w}_1 = \bar{M}^{-1}\bar{A}q_1$}}
    \State $\alpha_1 = \sqrt{\|w_{1,x}\|_{M}^2+\|w_{1,c}\|_{F^{-1}}^2}$
    {\hfill \textcolor{gray}{$\alpha_1 = \|\bar{w}_1\|_{\bar{M}}$}}
    \State $v_{1,x} = w_{1,x}/\alpha_1$; $v_{1,c} = w_{1,c}/\alpha_1$
    {\hfill \textcolor{gray}{$\bar{v}_1 = \bar{w}_1/\alpha_1$}}
    \State $k=1$
    \While{$k<$ maxit}
      \State $\hat{g}_k=N^{-1}(A^Tv_{k,x}+E^Tv_{k,c})$
      {\hfill \textcolor{gray}{$\hat{g}_k=N^{-1}\bar{A}^T\bar{v}_k$}}
      \State $h_{k}=Q_k^TN\hat{g}_k$; $g_k=\hat{g}_k-Q_kh_{k}$; $\beta_{k+1}=\|g_k\|_{N}$
      \State $q_{k+1}=g_k/\beta_{k+1}$; $Q_{k+1}=[Q_k,q_{k+1}]$
      \State $w_{k+1,x}=M^{-1}(Aq_{k+1}-\beta_{k+1}Mv_{k,x})$
      {\hfill \textcolor{gray}{$\bar{w}_{k+1}=\bar{M}^{-1}(\bar{A}q_{k+1}-\beta_{k+1}\bar{M}\bar{v}_k)$}}
      \State $w_{k+1,c}=F(E^Tq_{k+1}-\beta_{k+1}F^{-1}v_{k,c})$
      \State $\alpha_{k+1}=\sqrt{\|w_{k+1,x}\|_{M}^2+\|w_{k+1,c}\|_{F^{-1}}^2}$
      {\hfill \textcolor{gray}{$\alpha_{k+1}=\|\bar{w}_{k+1}\|_{\bar{M}}$}}
      \State $v_{k+1,x}=w_{k+1,x}/\alpha_{k+1}$; $v_{k+1,c}=w_{k+1,c}/\alpha_{k+1}$
      {\hfill \textcolor{gray}{$\bar{v}_{k+1}=\bar{w}_{k+1}/\alpha_{k+1}$}}
      \State $k=k+1$
    \EndWhile
  \end{algorithmic}
\end{algorithm}

From (\ref{e28}), the centered-preconditioned Schur complement $N^{-\frac{1}{2}}SN^{-\frac{1}{2}}$ with $S=\bar{A}^T\bar{M}^{-1}\bar{A}$ is reduced to an upper Hessenberg form as
\begin{equation}\label{e210}
\begin{array}{ll}
N^{-\frac{1}{2}}SN^{-\frac{1}{2}}N^{\frac{1}{2}}Q_k & =N^{-\frac{1}{2}}\bar{A}^T\bar{M}^{-1}\bar{A}Q_k =N^{-\frac{1}{2}}\bar{A}^T\bar{V}_kB_k \\
& =N^{-\frac{1}{2}}(NQ_kH_k+\beta_{k+1}Nq_{k+1}e_k^T)B_k \\
& =N^{\frac{1}{2}}Q_k(H_kB_k)+\alpha_k\beta_{k+1}N^{\frac{1}{2}}q_{k+1}e_k^T,
\end{array}
\end{equation}
where $H_kB_k$ is an upper Hessenberg matrix given by the Arnoldi process. 
When $k = n$, the decomposition residual $\beta_{k+1}Nq_{k+1}e_k^T$ in (\ref{e28}) vanishes and we obtain the upper Hessenberg form of $N^{-\frac{1}{2}}SN^{-\frac{1}{2}}=(N^{\frac{1}{2}}Q_n)(H_nB_n)(N^{\frac{1}{2}}Q_n)^T$. Therfore, we can implicitly obtain the upper Hessenberg form of $N^{-\frac{1}{2}}SN^{-\frac{1}{2}}$ 
in the sense that we find $Q_n$, $H_n$ and $B_n$ without knowledge of $S$.
Since $N^{-\frac{1}{2}}SN^{-\frac{1}{2}}$ is nonsymmetric, the column vectors in basis $Q_k$ generated in Arnoldi process 
satisfy the long recurrences, where one vector needs to be orthogonalized against all the previous vectors w.r.t. $\langle  \cdot, \cdot \rangle_{N}$. Thus, it is necessary to store all the vectors in $Q_k$ and use them in the orthogonalization process, to maintain global mutual orthogonality w.r.t. $\langle  \cdot,  \cdot \rangle_{N}$. 

Let $\bar{V}_{k}$ in (\ref{e28}) is defined as in (\ref{e26}). Then by (\ref{e17}) and (\ref{e26}), the (\ref{e28}) can be rewritten as
\begin{equation}\label{e211}
\left\{
\begin{array}{ll}
AQ_k=MV_{k,x}B_k,\quad EQ_k=F^{-1}V_{k,c}B_k, & V_{k,x}^TMV_{k,x}+V_{k,c}^TF^{-1}V_{k,c}=L_k, \\
A^TV_{k,x}+E^TV_{k,c}=NQ_kH_k+\beta_{k+1}Nq_{k+1}e_k^T, & Q_k^TNQ_k=I_k.
\end{array}
\right.
\end{equation}

In Section \ref{Sec4}, we will show how the decomposition
can be further reformulated by referring to the original system (\ref{e11}), thus avoiding the use of $E$ and $F$ and the factorization (\ref{e15}).

\section{The CRAIG algorithm for symmetric generalized saddle point problems (\ref{e11})}\label{Sec3}

\quad In this section, we firstly reformulate the GKB of the (1,2)-block of symmetric saddle point problem (\ref{e17}) given in (\ref{e23}) or  (\ref{e27}) as the GKB of the augmentation of the (1,2)-block of the  equivalent symmetric generalized saddle point problem (\ref{e11}).
Then, on this basis, we introduce the CRAIG algorithm for the symmetric generalized saddle point problem (\ref{e11}), including the corresponding linear solver, its stopping criteria and minimization property. Finally, we show that the proposed CRAIG algorithm is indeed theoretically equivalent to the SCR method where the preconditioned CG is applied to the related Schur-complement equations.

\subsection{The GKB for the symmetric generalized saddle point problem  (\ref{e11})}\label{SubSec31}


\quad We now reformulate Algorithm \ref{alg_S_GKB} without using $E$ and $F$. Define for $k=1,2,\ldots,$
\begin{equation}\label{e31}
\begin{array}{c}
w_{k}=w_{k,x}, \quad v_{k}=v_{k,x}, \quad
s_k=E^Tw_{k,c}, \quad t_k=E^Tv_{k,c}. \\
\end{array}
\end{equation}
Firstly, let
\begin{equation}\label{e32}
r_{1} = q_1,
\end{equation}
then line 2 of Algorithm \ref{alg_S_GKB}, (\ref{e15}) and (\ref{e31}) yield
\begin{equation}\label{e33}
\begin{array}{l}
w_{1,x} = w_1 = M^{-1}Aq_1, \quad w_{1,c} = FEq_1= FEr_1, \quad s_{1} = E^Tw_{1,c}=Cr_1.\\
\end{array}
\end{equation}
Line 3 of Algorithm \ref{alg_S_GKB}, (\ref{e31}) and (\ref{e33}) yield
\begin{equation}\label{e34}
\begin{array}{l}
\alpha_1 = \sqrt{\|w_{1,x}\|_{M}^2+\|w_{1,c}\|_{F^{-1}}^2}=\sqrt{\|w_1\|_{M}^2+r_1^TE^Tw_{1,c}}=\sqrt{\|w_1\|_{M}^2+r_1^Ts_1}.\\
\end{array}
\end{equation}
Line 4 of Algorithm \ref{alg_S_GKB}, (\ref{e15}) and (\ref{e31}) and (\ref{e33}) yield
\begin{equation}\label{e35}
\begin{array}{l}
v_{1,x}=v_{1} = w_{1}/\alpha_1, \quad v_{1,c} = w_{1,c}/\alpha_1=FEr_1/\alpha_1, 
\quad t_{1} = E^Tv_{1,c} = Cr_{1}/\alpha_1 = s_{1}/\alpha_1. 
\end{array}
\end{equation}
Lines 7 and 9 of Algorithm \ref{alg_S_GKB} and (\ref{e31}) yield
\begin{equation}\label{e36}
g_1=N^{-1}(A^Tv_{1}+t_{1}-\alpha_1Nq_1),\quad w_{2,x}=w_{2}=M^{-1}(Aq_{2}-\beta_{2}Mv_{1}).
\end{equation}
Define
\begin{equation}\label{e37}
r_{2}=q_{2}-\frac{\beta_{2}}{\alpha_1}r_{1}.
\end{equation}
then line 10 of Algorithm \ref{alg_S_GKB}, (\ref{e15}), (\ref{e31}) and (\ref{e35}) yield
\begin{equation}\label{e38}
w_{2,c}=FE\left(q_{2}-\frac{\beta_{2}}{\alpha_1}r_{1}\right)=FEr_{2}, \quad
s_{2} = E^Tw_{2,c} =Cr_{2}.
\end{equation}
Line 11 of Algorithm \ref{alg_S_GKB}, (\ref{e31}) and (\ref{e38}) yield
\begin{equation}\label{e39}
\begin{array}{c}
\alpha_{2}=\sqrt{\|w_{2,x}\|_{M}^2+\|w_{2,c}\|_{F^{-1}}^2}
=\sqrt{\|w_{2}\|_{M}^2+r_{2}^TE^Tw_{2,c}}=\sqrt{\|w_{2}\|_{M}^2+r_{2}^Ts_{2}}.\\
\end{array}
\end{equation}
Line 12 of Algorithm \ref{alg_S_GKB}, (\ref{e15}), (\ref{e31}) and (\ref{e38}) 
yield
\begin{equation}\label{e310}
\begin{array}{c}
v_{2,x}=v_{2}=w_{2}/\alpha_{2}, \quad v_{2,c}=w_{2,c}/\alpha_{2}=FEr_2/\alpha_2, 
\quad
t_{2}= E^Tv_{2,c}=Cr_{2}/\alpha_{2}=s_{2}/\alpha_{2}
. \\
\end{array}
\end{equation}

An induction argument shows that for all $k \geq 1$
\begin{equation}\label{e311}
\begin{array}{c}
g_k=N^{-1}(A^Tv_{k}+t_{k}-\alpha_kNq_k), \\
w_{k+1,x}=w_{k+1}=M^{-1}(Aq_{k+1}-\beta_{k+1}Mv_{k}). \\
\end{array}
\end{equation}
Define
\begin{equation}\label{e312}
r_{k+1}=q_{k+1}-\frac{\beta_{k+1}}{\alpha_k}r_{k},
\end{equation}
we obtain
\begin{equation}\label{e313}
\begin{array}{c}
w_{k+1,c}=FE(q_{k+1}-\frac{\beta_{k+1}}{\alpha_k}r_{k})=FEr_{k+1}, \quad
s_{k+1} = E^Tw_{k+1,c}=Cr_{k+1},\\
\alpha_{k+1}=\sqrt{\|w_{k+1,x}\|_{M}^2+\|w_{k+1,c}\|_{F^{-1}}^2}=\sqrt{\|w_{k+1}\|_{M}^2+r_{k+1}^TE^Tw_{k+1,c}}
=\sqrt{\|w_{k+1}\|_{M}^2+r_{k+1}^Ts_{k+1}}.\\
\end{array}
\end{equation}
Furthermore,
we obtain
\begin{equation}\label{e314}
\begin{array}{c}
v_{k+1,x}=v_{k+1}=w_{k+1}/\alpha_{k+1}, \quad v_{k+1,c}=w_{k+1,c}/\alpha_{k+1}=FEr_{k+1}/\alpha_{k+1}
, \\ t_{k+1}=E^Tv_{k+1,c}=Cr_{k+1}/\alpha_{k+1}=s_{k+1}/\alpha_{k+1}
. \\
\end{array}
\end{equation}
Thus, we obtain the steps 1-4, 8-9, 13-17 in Algorithm \ref{alg_CRAIG_G_SGSPP}.

Let
\begin{equation}\label{e315}
\begin{array}{c}
V_k=[v_{1},v_{2},\dots, v_{k}],\quad  
D_k=[\frac{1}{\alpha_1} r_{1},\frac{1}{\alpha_2} r_{2},\dots, \frac{1}{\alpha_k} r_{k}],
\quad T_k= [t_{1},t_{2},\dots,t_{k}].
\\
\end{array}
\end{equation}
This, along with  (\ref{e26}), (\ref{e31}), (\ref{e35}), (\ref{e310}) and (\ref{e314}), yields that
\begin{equation}\label{e316}
\begin{array}{c}
V_{k,x}=V_k,\quad
V_{k,c}=FED_k, \quad T_k=E^TV_{k,c}=CD_k. \\
\end{array}
\end{equation}
By (\ref{e32}), (\ref{e37}) and (\ref{e312}), we have
\begin{equation*}
\begin{array}{c}
q_{1}=r_1= (\frac{1}{\alpha_{1}}r_{1})\cdot \alpha_{1},\quad
q_{k+1}=(\frac{1}{\alpha_{k}}r_{k})\cdot\beta_{k+1}+ (\frac{1}{\alpha_{k+1}}r_{k+1})\cdot\alpha_{k+1}
,\quad k>0, \\
\end{array}
\end{equation*}
then combining with the structure of $B_k$ in (\ref{e24}) and the definition of $D_k$ in (\ref{e315}), we have
\begin{equation}\label{e317}
Q_k = D_kB_k.
\end{equation}
Then, by (\ref{e15}) and (\ref{e316}), the (\ref{e27}) can be rewritten as
\begin{equation}\label{e318}
\left\{
\begin{array}{ll}
AQ_k=MV_{k}B_k,\quad CQ_k=CD_kB_k=T_kB_k,\quad 
& V_{k}^TMV_{k}+D_{k}^TT_k=I_k,\\
A^TV_{k}+T_k=NQ_kB_k^T+\beta_{k+1}Nq_{k+1}e_k^T, & Q_k^TNQ_k=I_k,
\end{array}
\right.
\end{equation}
or equivalently,
\begin{equation}\label{e319}
\left\{
\begin{array}{ll}
\left[
\begin{array}{l}
A \\
C
\end{array}
\right]Q_k=\left[
\begin{array}{cc}
M& O \\
O & I_n
\end{array}
\right]\left[
\begin{array}{l}
V_k \\
T_k
\end{array}
\right]B_k, & [V_k^T\ D_{k}^T]\left[
\begin{array}{cc}
M& O \\
O& I_n
\end{array}
\right]\left[
\begin{array}{l}
V_k \\
T_k
\end{array}
\right]=I_k,\\
\left[A^T \ I_n\right] \left[
\begin{array}{l}
V_k \\
T_k
\end{array}
\right]=NQ_kB_k^T+\beta_{k+1}Nq_{k+1}e_k^T, & Q_k^TNQ_k=I_k,
\end{array}
\right.
\end{equation}
without explicitly including $V_{k,c}$. Hence,
the bidiagonalization (\ref{e318}) of the augmentation of (1,2)-block of the symmetric generalized saddle point problems (\ref{e11}) is the bidiagonalization (\ref{e319}) of the (1,2)- and (2,1)-blocks of the nonsymmetric saddle point problems (\ref{e18}). Consequently, the following algorithm based on the bidiagonalization (\ref{e319}) for (\ref{e18}) is the CRAIG method we want to propose for the symmetric generalized saddle point problems (\ref{e11}).


\subsection{The CRAIG algorithm for the generalized saddle point problem (\ref{e11})}\label{SubSec32}

\quad Using the $N$-orthogonality of the columns of $Q_{k+1}$,
the choice for $q_1$ given in (\ref{e21}) and the relations in (\ref{e319}), we can transform the system (\ref{e18}) into a simpler form
\begin{equation*}
\left[\begin{array}{cc}
		[V_k^T\ D_{k}^T] & \\
		& Q_k^T\\
	\end{array}\right]
\left[\begin{array}{cc|c}
		M && A\\
		&I_n& C\\
\hline
		A^T &I_n & \\
	\end{array}\right]\left[\begin{array}{cc}
		\left[
\begin{array}{l}
V_k \\
T_k
\end{array}
\right] & \\
		& Q_k\\
	\end{array}\right]\left[\begin{array}{c}
		z_k \\ y_k
	\end{array}\right]
	=\left[\begin{array}{cc}
		[V_k^T\ D_{k}^T] & \\
		& Q_k^T\\
	\end{array}\right]\left[\begin{array}{c}
		0 \\ 0 \\ b
	\end{array}\right],
\end{equation*}
or equivalently,
\begin{equation}\label{e320}
\left[\begin{array}{c c}
		I_k & B_k\\
		B_k^T & O\\
	\end{array}\right]
	\left[\begin{array}{c}
		z_k \\ y_k
	\end{array}\right]
=\left[\begin{array}{c}
		0 \\ \beta_1e_1
	\end{array}\right].
\end{equation}
The solution of (\ref{e320}) are then given by
\begin{equation}\label{e321}
z_k=\beta_1B_k^{-T}e_1,  \quad y_k= -B_k^{-1}z_k.
\end{equation}
We can build the $k$th approximate solution to (\ref{e18}) as
\begin{equation}\label{e322}
u^{(k)}=V_kz_k,\quad a^{(k)}=T_kz_k, \quad p^{(k)}=Q_ky_k.
\end{equation}
By (\ref{e317}), (\ref{e318}), (\ref{e321}) and (\ref{e322}), it is easy to obtain the relationship between $u^{(k)}$ and $p^{(k)}$, and that between $a^{(k)}$ and $p^{(k)}$:
\begin{equation}\label{e323}
\begin{array}{c}
u^{(k)}=V_kz_k=-V_kB_ky_k=-M^{-1}AQ_ky_k=-M^{-1}Ap^{(k)}, \\
a^{(k)}=T_kz_k=-T_kB_ky_k=-CQ_ky_k=-Cp^{(k)}.
\end{array}
\end{equation}

Given the structure of $\beta_1e_1$ and $B_k$, through (\ref{e321}), we find recursively
\begin{equation}\label{e324}
\zeta_1=\frac{\beta_1}{\alpha_1},\quad \zeta_{k}=-\frac{\beta_k}{\alpha_k}\zeta_{k-1},\quad 
z_k=\left[\begin{array}{c} z_{k-1}\\ \zeta_{k} \end{array}\right].
\end{equation}
Then, along with (\ref{e317}), the recursive update formulas for $u$, $a$ and $p$ are
\begin{equation}\label{e325}
\begin{array}{c}
u^{(k)}=V_kz_k=V_{k-1}z_{k-1}+v_{k}\zeta_{k}=u^{(k-1)}+\zeta_kv_{k}, \quad
a^{(k)}=T_kz_k=
a^{(k-1)}+\zeta_kt_{k}, \\
p^{(k)}=Q_ky_k= -Q_kB_k^{-1}z_k= -D_kz_k = -D_{k-1}z_{k-1}-(\frac{1}{\alpha_k} r_{k})\zeta_{k} = p^{(k-1)} - \frac{\zeta_k}{\alpha_k}r_k,
\end{array}
\end{equation}
with $u^{(0)}=0$, $a^{(0)}=0$ and $p^{(0)}=0$. In steps 5 and 18 in Algorithm \ref{alg_CRAIG_G_SGSPP}, only $u$ and $p$ need to be updated since the approximate solution of the original generalized saddle point problem (\ref{e11}) has already been obtained.

\begin{algorithm}[!ht]
\caption{The CRAIG algorithm for symmetric generalized saddle point problems (\ref{e11})}
\label{alg_CRAIG_G_SGSPP}
\textbf{Require}: 
$M\in R^{m\times m}$ SPD, $A\in R^{m\times n}$ with full column rank, $C\in R^{n\times n}$ SPSD, $b\in R^{n}$, maxit, tol
\\
\textbf{Output}: $u_{k}$, $p_{k}$
  \begin{algorithmic}[1]
    \State $\beta_1=\|b\|_{N^{-1}}$; $q_1=N^{-1}b/\beta_1$
    \State $w_{1} = M^{-1}Aq_1$; $r_1=q_1$; $s_1=Cr_1$
    {\hfill \textcolor{gray}{$w_{1,x} = M^{-1}Aq_1$; $w_{1,c} = FEq_1$}}
    \State $\alpha_1 = \sqrt{\|w_1\|_{M}^2+r_1^Ts_1}$
    {\hfill \textcolor{gray}{$\alpha_1 = \sqrt{\|w_{1,x}\|_{M}^2+\|w_{1,c}\|_{F^{-1}}^2}$}}
    \State $v_{1} = w_{1}/\alpha_1$; 
    $t_1=s_1/\alpha_1$
    {\hfill \textcolor{gray}{$v_{1,x} = w_{1,x}/\alpha_1$; $v_{1,c} = w_{1,c}/\alpha_1$}}
    \State $\zeta_1=\frac{\beta_1}{\alpha_1}$; $u^{(1)}=\zeta_1v_1$;
    $p^{(1)}=-\frac{\zeta_1}{\alpha_1}r_1$
    \State 
    $k=1$
    \While{
    $k<$ maxit}
      \State $g_k=N^{-1}(A^Tv_{k}+t_{k}-\alpha_kNq_k)$
      {\hfill \textcolor{gray}{$g_k=N^{-1}(A^Tv_{k,x}+E^Tv_{k,c}-\alpha_kNq_k)$}}
      \State $\beta_{k+1}=\|g_k\|_{N}$
     \If{$\frac{\beta_{k+1}}{\beta_{1}} |\zeta_k|<$ tol}{\hfill \textcolor{gray}{Stopping criterion}}
        \State \textbf{break};
     \EndIf
      \State $q_{k+1}=g_k/\beta_{k+1}$
      \State $w_{k+1}=M^{-1}(Aq_{k+1}-\beta_{k+1}Mv_{k})$
      {\hfill \textcolor{gray}{$w_{k+1,x}=M^{-1}(Aq_{k+1}-\beta_{k+1}Mv_{k,x})$}}
      \State $r_{k+1}=q_{k+1}-\frac{\beta_{k+1}}{\alpha_k}r_{k}$;
      $s_{k+1}=Cr_{k+1}$
      {\hfill \textcolor{gray}{ $w_{k+1,c}=F(E^Tq_{k+1}-\beta_{k+1}F^{-1}v_{k,c})$}}
      \State $\alpha_{k+1}=\sqrt{\|w_{k+1}\|_{M}^2+r_{k+1}^Ts_{k+1}}$
      {\hfill \textcolor{gray}{$\alpha_{k+1}=\sqrt{\|w_{k+1,x}\|_{M}^2+\|w_{k+1,c}\|_{F^{-1}}^2}$}}
      \State $v_{k+1}=w_{k+1}/\alpha_{k+1}$; 
      $t_{k+1}=s_{k+1}/\alpha_{k+1}$
      \State $\zeta_{k+1}=-\frac{\beta_{k+1}}{\alpha_{k+1}}\zeta_{k}$; $u^{(k+1)}=u^{(k)}+\zeta_{k+1}v_{k+1}$; $p^{(k+1)}=p^{(k)}-\frac{\zeta_{k+1}}{\alpha_{k+1}}r_{k+1}$
      \State $k=k+1$
    \EndWhile
  \end{algorithmic}
\end{algorithm}

Note that Algorithm \ref{alg_CRAIG_G_SGSPP} does not contain references to $E$ and $F$ and it works directly with the formulation (\ref{e11}). This is its main advantage.
When $C=0$, Algorithm \ref{alg_CRAIG_G_SGSPP} reduces to 
the CRAIG solver for symmetric saddle point systems given in \cite{MA2013,VDADCKUR2023}.
In each step of Algorithm \ref{alg_CRAIG_G_SGSPP}, 
one more matrix-vector product $s_k=Cr_k$ and one more scalar product $r_{k}^Ts_{k}$ are required than CRAIG. Moreover, Algorithm \ref{alg_CRAIG_G_SGSPP} requires two more vectors of storage: $s$ and $t$ than CRAIG. When $n$ is much smaller than $m$, the amount of computation and storage space introduced in Algorithm \ref{alg_CRAIG_G_SGSPP} than CRAIG is limited.
In the algorithm, it is sufficient to store only the latest left vectors $v_{k}$ and $t_{k}$ to compute $v_{k+1}$ and $t_{k+1}$, and store only the latest right vectors $q_{k}$ to compute $q_{k+1}$, similar to the case of $C=O$.

Next, we give error estimates for the errors on $u-u^{(k)}$ and $p-p^{(k)}$, and on the dual norm of the residual $r^{(k)} =b - A^T u^{(k)} + Cp^{(k)}$.

By (\ref{e316}), (\ref{e317}), (\ref{e318}), (\ref{e321}) and (\ref{e322}), at step $k$ of Algorithm \ref{alg_CRAIG_G_SGSPP}, we have
\begin{equation}\label{e326}
\begin{array}{l}
\|u-u^{(k)}\|_{M}^2+(p-p^{(k)})^TC(p-p^{(k)}) \\
\quad =\|V_nz_n-V_kz_k\|_{M}^2+(Q_ny_n-Q_ky_k)^TC(Q_ny_n-Q_ky_k)\\
\quad =\|V_nz_n-V_kz_k\|_{M}^2+(Q_nB_n^{-1}z_n-Q_kB_k^{-1}z_k)^TC(Q_nB_n^{-1}z_n-Q_kB_k^{-1}z_k)\\
\quad =\|V_nz_n-V_kz_k\|_{M}^2+(D_nz_n-D_kz_k)^TC(D_nz_n-D_kz_k)\\
\quad =\left(z_n-\left[\begin{array}{c} z_k\\ O_{n-k,1} \end{array}\right]\right)^T(V_n^TMV_n+D_n^TCD_n)\left(z_n-\left[\begin{array}{c} z_k\\ O_{n-k,1} \end{array}\right]\right)=[O\ z_{n-k}^T]\left[\begin{array}{c} O\\ z_{n-k} \end{array}\right]\\
\quad =z_{n-k}^Tz_{n-k}\\
\quad =\sum_{i=k+1}^n\zeta_i^2,
\end{array}
\end{equation}
where $V_n=[V_{k}\ V_{n-k}]\in R^{m\times n}$ with $V_{n-k}
\in R^{m\times (n-k)}$, $D_n=[D_{k}\ D_{n-k}]\in R^{n\times n}$ with $D_{n-k}
\in R^{n\times (n-k)}$,  and $z_n=[z_{k}; z_{n-k}]\in R^{n}$ with $z_{n-k}=[\zeta_{k+1},\zeta_{k+2},\cdots,\zeta_{n}]^T\in R^{n-k}$.
It is the error in the energy norm for the primal variable $\bar{u}$ of (\ref{e17}), i.e., $\|\bar{u}-\bar{u}_k\|_{\bar{M}}$ rather than that for the primal variable $u$ of (\ref{e11}).
In addition, by (\ref{e323}), we have
\begin{equation}\label{e327}
\begin{array}{ll}
\|u-u^{(k)}\|_{M}^2+(p-p^{(k)})^TC(p-p^{(k)})&=\|M^{-1}A(p-p^{(k)})\|_{M}^2+(p-p^{(k)})^TC(p-p^{(k)})\\
&=\|p-p^{(k)}\|_{A^TM^{-1}A}^2+(p-p^{(k)})^TC(p-p^{(k)})\\
&=\|p-p^{(k)}\|_{A^TM^{-1}A+C}^2=\|p-p^{(k)}\|_{S}^2. \\
\end{array}
\end{equation}
Since the matrix $C$ is symmetric positive semidefinite, by (\ref{e326}) and (\ref{e327}), we have
\begin{equation}\label{e328}
\|u-u^{(k)}\|_{M}^2\leq \|p-p^{(k)}\|_{A^TM^{-1}A+C}^2 =
\sum_{i=k+1}^n\zeta_i^2.
\end{equation}

If we truncate the sum above to only its first $d$ terms, we get a lower bound on (the upper bound of) the energy norm of the error for $u$ and $p$. The subscript $d$ stands for delay, because we can compute this lower bound corresponding to a given step $k$ only after an additional $d$ steps
\begin{equation*}
\xi_{k,d}^2=\sum_{i=k+1}^{k+d}\zeta_i^2<\sum_{i=k+1}^n\zeta_i^2.
\end{equation*}
With this bound for the absolute error, we can devise one for the relative error in (\ref{e329}), which can be used as a stopping criterion of Algorithm \ref{alg_CRAIG_G_SGSPP}
\begin{equation}\label{e329}
\bar{\xi}_{k,d}^2=\frac{\sum_{i=k-d+1}^{k}\zeta_i^2}{\sum_{i=1}^{k}\zeta_i^2}.
\end{equation}

By (\ref{e21}), (\ref{e318}), (\ref{e321}), (\ref{e322}) and (\ref{e324}), we have
\begin{equation}\label{e330}
\begin{array}{ll}
r^{(k)}_{CRAIG} =b-A^Tu^{(k)} +Cp^{(k)} & =\beta_1Nq_1-A^TV_kz_k+CQ_ky_k  \\ &=\beta_1Nq_1-A^TV_kz_k+T_kB_ky_k \\
& = \beta_1Nq_1-(A^TV_k+T_k)z_k \\
&= \beta_1Nq_1-(NQ_kB_k^T+\beta_{k+1}Nq_{k+1}e_k^T)z_k \\
&= \beta_1NQ_ke_1 -\beta_1NQ_ke_1- \beta_{k+1}Nq_{k+1}e_k^Tz_k \\
&= - \zeta_k\beta_{k+1}Nq_{k+1}.  \\
\end{array}
\end{equation}
It follows from (\ref{e330}) that the residual of the second equation of (\ref{e11}) is parallel to the vector $Nq_{k+1}$ and it is orthogonal to all right vectors stored in $Q_k$, i.e., $(r^{(k)})^TQ_k=0$ owing to the $N$-orthogonality of $q_{k+1}$.
Combining the first equation in (\ref{e323}) and (\ref{e330}), we have the residual of (\ref{e11}) at $(u^{(k)},p^{(k)})$:
\begin{equation}\label{e331}
\mathbf{res}_k^{CRAIG} \triangleq
\left[\begin{array}{c}
0\\ b
\end{array}\right]-\left[\begin{array}{cc}
M& A\\
A^T & -C\\
\end{array}\right]\left[\begin{array}{c}
u^{(k)} \\ p^{(k)}
\end{array}\right]=
\left[\begin{array}{c}
O \\ r^{(k)}_{CRAIG}
\end{array}\right]
=\left[\begin{array}{c}
0\\ - \zeta_k\beta_{k+1}Nq_{k+1}
\end{array}\right].
\end{equation}
Then,
at step $k$ of Algorithm  \ref{alg_CRAIG_G_SGSPP}, we have the following residual norm
\begin{equation}\label{e332}
\|\mathbf{res}_k^{CRAIG}\|_{D_0^{-1}}=\|r^{(k)}_{CRAIG}\|_{N^{-1}}=
\|b-A^Tu^{(k)}+Cp^{(k)}\|_{N^{-1}}=\beta_{k+1}|\zeta_k|,
\end{equation}
where $D_0=\mbox{blkdiag}(M,N)$ and the $N$-orthogonality of $q_{k+1}$ is used. The residual norm (\ref{e332}) in the case of $C=O$ was given in the equation (4.3) of \cite{MA2013}.
With this bound for the absolute residual, we can devise one for the relative residual in (\ref{e333}), which can be used as another stopping criterion for Algorithm \ref{alg_CRAIG_G_SGSPP}
\begin{equation}\label{e333}
\frac{\|\mathbf{res}_k^{CRAIG}\|_{D_0^{-1}}}{\|\mathbf{res}_0^{CRAIG}\|_{D_0^{-1}}} =
\frac{\|r^{(k)}_{CRAIG}\|_{N^{-1}}}{\|r^{(0)}_{CRAIG}\|_{N^{-1}}}
=\frac{\|b-A^Tu^{(k)}+Cp^{(k)}\|_{N^{-1}}}{\|b\|_{N^{-1}}}=\frac{\beta_{k+1}}{\beta_{1}}|\zeta_k|.
\end{equation}
Since $\beta_{k+1}$ and $\zeta_k$ can be recursively calculated, the relative residual norm can be computed each step very cheaply. In specific numerical experiments, for simplicity, we use (\ref{e333}) as the stopping rule, see steps 10-12 in Algorithm \ref{alg_CRAIG_G_SGSPP}. It is worth emphasizing that although using (\ref{e333}) is more attractive than using (\ref{e329}), unlike the energy norm of the error, the residual dual norm does not have any physical meaning and the residual dual norm may be not a reliable indicator of the error if the matrix condition number is large \cite{ECJLZS2024}.

\begin{remark} \label{R31}
We observe that, owing to the nonsingularity of both $M$ and $N$ and
the full rank of $A$, all $\beta_k\geq 0$ and $\alpha_k>0$, $(k = 1,...,n)$. 
Since $M\succ O$, $C\succeq O$, and $A$ has full column rank, the linear system (\ref{e11}) has only one solution and $b\in \mbox{Range}([A^T\ C])$. Assume that $\alpha_{k+1}=0$, $k\in \{0,1,\cdots,n-1\}$, then (\ref{e319}) becomes
\begin{equation}\label{e334}
\left\{
\begin{array}{ll}
\left[
\begin{array}{l}
A \\
C
\end{array}
\right]Q_{k+1}=\left[
\begin{array}{cc}
M& O \\
O & I_n
\end{array}
\right]\left[
\begin{array}{l}
V_k \\
T_k
\end{array}
\right][B_k\ \beta_{k+1}e_k], & [V_k^T\ D_{k}^T]\left[
\begin{array}{cc}
M& O \\
O& I_n
\end{array}
\right]\left[
\begin{array}{l}
V_k \\
T_k
\end{array}
\right]=I_k,\\
\left[A^T \ I_n\right] \left[
\begin{array}{l}
V_k \\
T_k
\end{array}
\right]=NQ_{k+1}\left[
\begin{array}{l}
B_k^T \\
\beta_{k+1}e_k^T
\end{array}\right], & Q_{k+1}^TNQ_{k+1}=I_{k+1}.
\end{array}
\right.
\end{equation}
Let $[B_k\ \beta_{k+1}e_k]= \hat{U}_{k}[\Sigma_{k}\ O]\hat{W}_{k+1}^T$ be the singular value decomposition of $[B_k\ \beta_{k+1}e_k]$, where $\hat{U}_{k}\in R^{k\times k}$ is orthogonal, $\Sigma_{k}=\mbox{diag}(\sigma_1,\sigma_2,\dots, \sigma_k)\in R^{k\times k}$ with $\sigma_1\geq\sigma_2\geq\dots\geq\sigma_k>0$, $\hat{W}_{k+1}\in R^{(k + 1)\times (k + 1)}$ is orthogonal.
We have from (\ref{e334}) that
\begin{equation}\label{e335}
\left\{
\begin{array}{ll}
\left[
\begin{array}{l}
A \\
C
\end{array}
\right]\hat{Q}_{k+1}=\left[
\begin{array}{cc}
M& O \\
O & I_n
\end{array}
\right]\left[
\begin{array}{l}
\hat{V}_k \\
\hat{T}_k
\end{array}
\right][\Sigma_{k}\ O], & [\hat{V}_k^T\ \hat{D}_{k}^T]\left[
\begin{array}{cc}
M& O \\
O& I_n
\end{array}
\right]\left[
\begin{array}{l}
\hat{V}_k \\
\hat{T}_k
\end{array}
\right]=I_k,\\
\left[A^T \ I_n\right] \left[
\begin{array}{l}
\hat{V}_k \\
\hat{T}_k
\end{array}
\right]=N\hat{Q}_{k+1}\left[
\begin{array}{l}
\Sigma_{k} \\
O
\end{array}
\right], & \hat{Q}_{k+1}^TN\hat{Q}_{k+1}=I_{k+1},
\end{array}
\right.
\end{equation}
where
\begin{equation*}
\hat{Q}_{k+1}=Q_{k+1}\hat{W}_{k+1},\quad \hat{V}_k=V_k\hat{U}_{k},\quad \hat{T}_k=T_k\hat{U}_{k},\quad \hat{D}_k=D_k\hat{U}_{k}.
\end{equation*}
The relations (\ref{e335}) show that $\hat{q}_{k+1}=Q_{k+1}\hat{w}_{k+1} \in \mbox{Null}(\left[A^T \ C\right]^T)$ where $\hat{q}_{k+1}$ and $\hat{w}_{k+1}$ are the last columns of $\hat{Q}_{k+1}$ and $\hat{W}_{k+1}$, respectively.
Since $b\in \mbox{Range}([A^T\ C])$, then $q_1\in \mbox{Range}(N^{-1}[A^T\ C])=\mbox{Range}([A^T\ C])$. A recursion argument easily establishes that each $q_j\in \mbox{Range}([A^T\ C])$. In this case, the vector $\hat{q}_{k+1}=Q_{k+1}\hat{w}_{k+1}$, a combination of $q_j$, $j=1,2,\cdots, k+1$, thus lies in $\mbox{Range}([A^T\ C])$ which is in contradiction with the previous conclusion that $\hat{q}_{k+1}=Q_{k+1}\hat{w}_{k+1} \in \mbox{Null}(\left[A^T \ C\right]^T)$. Therefore, if $b\in \mbox{Range}([A^T\ C])$, then Algorithm \ref{alg_CRAIG_G_SGSPP} cannot terminate with $\alpha_{k+1} = 0$.
On the other hand, if $\beta_{k+1}=0$, then by (\ref{e331}), Algorithm \ref{alg_CRAIG_G_SGSPP} terminates with an exact solution.
An example of such early termination is if $b$ is an eigenvector of the right-preconditioned Schur complement $SN^{-1}$ corresponding to the eigenvalue $\alpha_1^2$, we have $\beta_1 =\|b\|_{N^{-1}}$, $q_1 = N^{-1}b/\beta_1 $, $q_2 = 0$, $\beta_2 = 0$ and the bidiagonalization terminates with the first iterates being the exact solution. The expression of the upper bound of the $M$-norm of the error on $u$ and the $S$-norm of the error on $p$ in (\ref{e328}) entail that the sequence $\|p-p^{(k)}\|_{S}$ decreases strictly.
\end{remark}

From (\ref{e318}) or (\ref{e319}), we can draw a link between the bidiagonalization of the augmentation of the off-diagonal block $A$ of (\ref{e11}) or the bidiagonalization of the off-diagonal blocks $[A^T\ C]^T$ and $[A^T\  I_n]$ of (\ref{e18}) and the tridiagonalization of the centered-preconditioned Schur complement $N^{-\frac{1}{2}}SN^{-\frac{1}{2}}$ with $S=A^TM^{-1}A+C$, which is exactly the same as that given in (\ref{e25}).
The equation (\ref{e25}) is the matrix form of the Lanczos process in CG with zero initial guess for the centered-preconditioned system of the Schur-complement equation (\ref{e13}): 
\begin{equation}\label{e336}
N^{-\frac{1}{2}}SN^{-\frac{1}{2}}\hat{p}=-N^{-\frac{1}{2}}b, \quad p=N^{-\frac{1}{2}}\hat{p}.
\end{equation}

An important property of the Lanczos vectors in $N^{\frac{1}{2}}Q_k$ is that they lie in the Krylov subspace
\begin{equation}\label{e337}
\mathcal{K}_k\triangleq\mathcal{K}_k(N^{-\frac{1}{2}}SN^{-\frac{1}{2}},N^{-\frac{1}{2}}b)=\mbox{span}\{N^{-\frac{1}{2}}b,N^{-\frac{1}{2}}SN^{-\frac{1}{2}}N^{-\frac{1}{2}}b,\cdots,(N^{-\frac{1}{2}}SN^{-\frac{1}{2}})^{k-1}N^{-\frac{1}{2}}b\}.
\end{equation}
At iteration $k$ of the CG method for (\ref{e336}), we look for an approximate solution $\hat{p}^{(k)}_{CG}=N^{\frac{1}{2}}Q_ky_k$ characterized variationally via
\begin{equation}\label{e338}
\hat{p}^{(k)}_{CG}\in \mathcal{K}_k(N^{-\frac{1}{2}}SN^{-\frac{1}{2}},N^{-\frac{1}{2}}b), \quad \hat{r}^{(k)}_{CG}=-N^{-\frac{1}{2}}b-N^{-\frac{1}{2}}SN^{-\frac{1}{2}}\hat{p}^{(k)}_{CG}
\perp \mathcal{K}_k(N^{-\frac{1}{2}}SN^{-\frac{1}{2}},N^{-\frac{1}{2}}b).
\end{equation}
This implies that
$y_k$ is chosen to minimize the energy norm of the error within the Krylov subspace $\mathcal{K}_k$:
\begin{equation}\label{e339}
\|\hat{p}-\hat{p}_{CG}^{(k)}\|_{N^{-\frac{1}{2}}SN^{-\frac{1}{2}}}
=\min_{y\in R^k}
\|\hat{p}-N^{\frac{1}{2}}Q_ky\|_{N^{-\frac{1}{2}}SN^{-\frac{1}{2}}}=\|\hat{p}-N^{\frac{1}{2}}Q_ky_k\|_{N^{-\frac{1}{2}}SN^{-\frac{1}{2}}},
\end{equation}
where  
$y_k$ satisfies that
\begin{equation*}
B_k^TB_ky_k=-\beta_1e_1,
\end{equation*}
with $\beta_1=\|N^{-\frac{1}{2}}b\|_2=\|b\|_{N^{-1}}$. Thus, by the above analysis, we have the solution of (\ref{e13}) 
at the $k$th step of the CG with preconditioner $N$:
\begin{equation*}
p^{(k)}_{CG}=N^{-\frac{1}{2}}\hat{p}^{(k)}_{CG}=N^{-\frac{1}{2}}N^{\frac{1}{2}}Q_ky_k=-Q_k(B_k^{-1}(B_k^{-T}\beta_1e_1)),
\end{equation*}
which is exactly the lower block $p^{(k)}$ of the approximate solution to (\ref{e11}) given in (\ref{e322}). Then, by the first equation of (\ref{e323}),  we have the following result.

\begin{theorem}\label{T31}
	The CRAIG iterates $u^{(k)}$ produced in Algorithm \ref{alg_CRAIG_G_SGSPP} are related to the iterates $p^{(k)}$ of the conjugate gradient method applied to the Schur-complement equation (\ref{e13})
with preconditioner $N$ according to $u^{(k)}=-M^{-1}Ap^{(k)}$.
\end{theorem}

\begin{remark} \label{R32}
Note that when the (2,2)-block $C$ of the symmetric generalized saddle point problem (\ref{e11}) is also SPD, the class of systems 
are called symmetric quasi-definite systems. 
In \cite[Chapter 5]{DOMA2017}, a generalized CRAIG solver is proposed for the symmetric quasi-definite systems of the form (\ref{e11}), 
and it is theoretically equivalent to the SCR method with inner CG iteration applied to the related Schur-complement equations (\ref{e13}) with the SPD preconditioner $C$. Obviously, Algorithm \ref{alg_CRAIG_G_SGSPP} can be applied to the class of symmetric quasi-definite systems, but the method proposed in \cite{DOMA2017} can not be used to solve the generalized saddle point problem (\ref{e11}) if $C$ is only SPSD. Thus, one advantage of our proposed method over the generalized CRAIG method introduced in \cite{DOMA2017} is that the former has a wider range of applications than the latter.
Another advantage is that our proposed method can be more flexible in selecting the SPD preconditioner $N$, making it more effective than the method with the fixed SPD preconditioner $C$ proposed in \cite{DOMA2017}.
\end{remark}

Let $\mathcal{V}_k = \mbox{span}\{v_1, \dots, v_k\}$ and $\mathcal{Q}_k = \mbox{span}\{q_1, \dots, q_k\}$. Then, by (\ref{e322}), (\ref{e327}), (\ref{e330}), (\ref{e339}) and Theorem \ref{T31}, 
for any step $k$, we have that
\begin{equation}\label{e340}
\min_{\small{\begin{array}{c}
u^{(k)}\in \mathcal{V}_{k}, p^{(k)}\in \mathcal{Q}_{k}\\
\left(A^Tu^{(k)}-Cp^{(k)}-b\right)\perp \mathcal{Q}_k\\
\end{array}}
} \sqrt{\|u-u^{(k)}\|_{M}^2 
+\left(p-p^{(k)}\right)^TC\left(p-p^{(k)}\right)}
\end{equation}
is met for $u^{(k)}$ and $p^{(k)}$ as computed by Algorithm \ref{alg_CRAIG_G_SGSPP}. 
The error minimization property (\ref{e340}) for Algorithm \ref{alg_CRAIG_G_SGSPP} is similar to that for CRAIG proposed in \cite{MA2013}.

\begin{remark} \label{R33}
Theorem \ref{T31} shows that for the symmetric generalized saddle point problem (\ref{e11}), our proposed CRAIG and SCR(CG) produce the same iterates $p^{(k)}$ at each step, and also produce the same iterates $u^{(k)}$ when the termination of the iterations.
The fact is similar to the case of $C=O$ given in
\cite[Chapter 5]{DOMA2017} and \cite{MA2013}.
However, the proposed CRAIG computes iteratively the upper block $u^{(k)}$ of the approximation solution of (\ref{e11}), while the SCR(CG) computes $u^{(k)}$ by multiplying $A$ and $p^{(k)}$ and solving the system with $M$ by matrix factorization until the stopping condition has been fulfilled. This lead to that the total computational cost for our CRAIG may be smaller than the SCR(CG) when the scale of the system (\ref{e11}) is very large and the number of iteration of the CRAIG is 
very small, see the ``time" in the Tables \ref{t51}-\ref{t53}.
More importantly, the proposed CRAIG is an example of solvers originating from GKB like CRAIG, LSQR and LSMR \cite{MA2013,CCPMAS1982,DCLFMS2011}, thus its use can be more desirable than the alternative of using SCR(CG) with inner CG iteration directly on the Schur-complement equation, which would involve operating with a squared condition number and the difficulty in ensuring the accuracy of its solution, as described in \cite[Chapter 8.1]{YS2003}, and \cite{CCPMAS1982}.
\end{remark}

\section{The nsCRAIG algorithm for nonsymmetric generalized saddle point problems (\ref{e11})}\label{Sec4}

\quad In this section, We turn to the case where the (1,1)-block $M$ in (\ref{e11}) is not symmetric. Firstly, we reformulate the decomposition of the (1,2)-block of nonsymmetric saddle point problem (\ref{e17}) given in (\ref{e28}) or (\ref{e211}) as the decomposition of the augmentation of the (1,2)-block of the equivalent nonsymmetric generalized saddle point problems (\ref{e11}).
Then, on this basis, we introduce the nsCRAIG algorithm for the nonsymmetric generalized saddle point problems (\ref{e11}), including the corresponding linear solver and its stopping criteria.
Finally, we show that the proposed nsCRAIG algorithm is indeed theoretically equivalent to the SCR method where the preconditioned FOM is applied to the associated Schur-complement equations.


\subsection{The GKB for the nonsymmetric generalized saddle point problems (\ref{e11})}\label{SubSec41}

\quad Similar to the process given in Section \ref{SubSec31}, we can reformulate Algorithm \ref{alg_NS_GKB} by avoiding use of $E$ and $F$, 
and then list it in the steps 1-4, 7-8, 13-17 in Algorithm \ref{alg_nsCRAIG_G_NSGSPP}. 


Let $V_{k,x}$  and $V_{k,c}$ in (\ref{e211}) are defined as in (\ref{e316}). Then, by (\ref{e15}) and (\ref{e316}), the (\ref{e211}) can be rewritten as
\begin{equation}\label{e41}
\left\{
\begin{array}{ll}
AQ_k=MV_{k}B_k,\quad CQ_k=CD_kB_k=T_kB_k, & V_{k}^TMV_{k}+D_k^TT_k=L_k, \\
A^TV_{k}+T_k=NQ_kH_k+\beta_{k+1}Nq_{k+1}e_k^T, & Q_k^TNQ_k=I_k,
\end{array}
\right.
\end{equation}
or equivalently,
\begin{equation}\label{e42}
\left\{
\begin{array}{ll}
\left[
\begin{array}{l}
A \\
C
\end{array}
\right]Q_k=\left[
\begin{array}{cc}
M& O \\
O & I_n
\end{array}
\right]\left[
\begin{array}{l}
V_k \\
T_k
\end{array}
\right]B_k, & [V_k^T\ D_{k}^T]\left[
\begin{array}{cc}
M& O \\
O& I_n
\end{array}
\right]\left[
\begin{array}{l}
V_k \\
T_k
\end{array}
\right]=L_k,\\
\left[A^T \ I_n\right] \left[
\begin{array}{l}
V_k \\
T_k
\end{array}
\right]=NQ_kH_k+\beta_{k+1}Nq_{k+1}e_k^T, & Q_k^TNQ_k=I_k,
\end{array}
\right.
\end{equation}
without explicitly including $V_{k,c}$. Hence, the decomposition (\ref{e41}) of the augmentation of (1,2)-block of the nonsymmetric generalized saddle point problems (\ref{e11}) is the decomposition (\ref{e42}) of the (1,2)- and (2,1)- blocks of the nonsymmetric saddle point problems (\ref{e18}). Consequently, the following algorithm based on the decomposition (\ref{e42}) for (\ref{e18}) is the nsCRAIG solver we want to propose for the nonsymmetric generalized saddle point problems (\ref{e11}).

\subsection{The nsCRAIG algorithm for the generalized saddle point problems (\ref{e11})}\label{SubSec42}

\quad Using the $N$-orthogonality of the columns of $Q_{k+1}$, 
and the choice for $q_1$ given in (\ref{e21}), the relations in (\ref{e42}) and the relation $H_k=B_k^TL_k^T$, we can transform the system (\ref{e18}) into a simpler form
\begin{equation}\label{e43}
\left[\begin{array}{c c}
		L_k & L_kB_k\\
		H_k & O\\
	\end{array}\right]
	\left[\begin{array}{c}
		z_k \\ y_k
	\end{array}\right]=\left[\begin{array}{c c}
		L_k & L_kB_k\\
		B_k^TL_k^T & O\\
	\end{array}\right]
	\left[\begin{array}{c}
		z_k \\ y_k
	\end{array}\right]
	=\left[\begin{array}{c}
		0 \\ \beta_1e_1 
	\end{array}\right].
\end{equation}
The solution of (\ref{e43}) are then given by
\begin{equation}\label{e44}
z_k=\beta_1H_k^{-1}e_1=\beta_1L_k^{-T}B_k^{-T}e_1, \quad y_k= -B_k^{-1}z_k.
\end{equation}
We can build the $k$th approximate solution to (\ref{e18}) as
\begin{equation}\label{e45}
u^{(k)}=V_kz_k, \quad a^{(k)}=T_kz_k, \quad p^{(k)}=Q_ky_k,
\end{equation}
which is identical in form to (\ref{e322}) for symmetric case. By (\ref{e41}), (\ref{e44}) and (\ref{e45}), we have
\begin{equation}\label{e46}
\begin{array}{c}
u^{(k)}=-M^{-1}Ap^{(k)}, \quad a^{(k)}=-Cp^{(k)},
\end{array}
\end{equation}
which are exactly the same as the relations in (\ref{e323}) for symmetric case.

For exact arithmetic, the nsCRAIG algorithm gives the exact solution after at most $n$ steps.
In the absence of exact arithmetic and/or if we only have a partial decomposition ($k < n$), $u^{(k)}$, $a^{(k)}$ and $p^{(k)}$ are only approximate solutions.
The decomposition (\ref{e41}) or (\ref{e42}) progresses by relying on the vectors $v$, $t$ and $q$, without referring to the current iterates $u^{(k)}$, $a^{(k)}$ and $p^{(k)}$. Consequently, forming this approximation can be postponed until the stopping condition is fulfilled.
Specifically, once the stopping condition is met, by (\ref{e44}), 
we firstly only need to apply the inverses $B_k^{-T}$ and $L_k^{-T}$ once to obtain $z_k$ and then apply the inverse $B_k^{-1}$ once to obtain $y_k$. Owing to the particular structure of the matrices $B_k$ and $L_k$,
the first inverse requires forward substitution, while the last two inverses require backward substitution.
Secondly, the approximate solution $p^{(k)}$ to (\ref{e18}) or (\ref{e11}) is found by the last equation in (\ref{e45}) and then the approximate solutions $u^{(k)}$ and $a^{(k)}$ to (\ref{e18}) are yielded by (\ref{e46}) since $Q_k$ is stored but $V_k$ and $T_k$ are not stored. 
In step 10 of Algorithm \ref{alg_nsCRAIG_G_NSGSPP}, we only need to update $u$ and $p$, which has already been approximate solution of the original generalized saddle point problem (\ref{e11}).

\begin{algorithm}[!ht]
\caption{The nsCRAIG algorithm for the nonsymmetric generalized saddle point problems (\ref{e11})}
\label{alg_nsCRAIG_G_NSGSPP}
\textbf{Require}: 
$M\in R^{m\times m}$ NSPD, $A\in R^{m\times n}$ with full column rank, $C\in R^{n\times n}$ SPSD, $b\in R^{n}$, maxit, tol
\\
\textbf{Output}: $u^{(k)}$, $p^{(k)}$
  \begin{algorithmic}[1]
    \State $\beta_1=\|b\|_{N^{-1}}$; $q_1=N^{-1}b/\beta_1$
    \State $w_{1} = M^{-1}Aq_1$; $r_1=q_1$; $s_1=Cr_1$;
    {\hfill \textcolor{gray}{$w_{1,x} = M^{-1}Aq_1$; $w_{1,c} = FEq_1$}}
    \State $\alpha_1 = \sqrt{\|w_1\|_{M}^2+r_1^Ts_1}$
    {\hfill \textcolor{gray}{$\alpha_1 = \sqrt{\|w_{1,x}\|_{M}^2+\|w_{1,c}\|_{F^{-1}}^2}$}}
    \State $v_{1} = w_{1}/\alpha_1$; 
    $t_1=s_1/\alpha_1$
    {\hfill \textcolor{gray}{$v_{1,x} = w_{1,x}/\alpha_1$; $v_{1,c} = w_{1,c}/\alpha_1$}}
    \State $\chi_1=\frac{\beta_1}{\alpha_1}$; 
    $k=1$
    \While{$k<$ maxit}
      \State $\hat{g}_k=N^{-1}(A^Tv_{k}+t_{k})$
      {\hfill \textcolor{gray}{$\hat{g}_k=N^{-1}(A^Tv_{k,x}+E^Tv_{k,c})$}}
     \State $h_{k}=Q_k^TN\hat{g}_k$; $g_k=\hat{g}_k-Q_kh_{k}$; $\beta_{k+1}=\|g_k\|_{N}$
     \If{$\frac{\beta_{k+1}}{\beta_{1}} |\chi_k|<$ tol} {\hfill \textcolor{gray}{Stopping criterion}}
     \State $y_k=-B_k^{-1}(H_k^{-1}(\beta_1e_1))$; $p^{(k)}=Q_{k}y_k$; $u^{(k)}=-M^{-1}A p_{k}$
        \State \textbf{break}; 
     \EndIf
      \State $q_{k+1}=g_k/\beta_{k+1}$; $Q_{k+1}=[Q_k,q_{k+1}]$
      \State $w_{k+1}=M^{-1}(Aq_{k+1}-\beta_{k+1}Mv_{k})$
      {\hfill \textcolor{gray}{$w_{k+1,x}=M^{-1}(Aq_{k+1}-\beta_{k+1}Mv_{k,x})$;}}
      \State $r_{k+1}=q_{k+1}-\frac{\beta_{k+1}}{\alpha_k}r_{k}$; $s_{k+1}=Cr_{k+1}$
      {\hfill \textcolor{gray}{ $w_{k+1,c}=F(Eq_{k+1}-\beta_{k+1}F^{-1}v_{k,c})$;}}
      \State $\alpha_{k+1}=\sqrt{\|w_{k+1}\|_{M}^2+r_{k+1}^Ts_{k+1}}$
      {\hfill \textcolor{gray}{$\alpha_{k+1}=\sqrt{\|w_{k+1,x}\|_{M}^2+\|w_{k+1,c}\|_{F^{-1}}^2}$}}
      \State $v_{k+1}=w_{k+1}/\alpha_{k+1}$; 
      $t_{k+1}=s_{k+1}/\alpha_{k+1}$
      \State $\chi_{k+1}=-\frac{\beta_{k+1}}{\alpha_{k+1}}\chi_{k}$
      \State $k=k+1$
    \EndWhile
  \end{algorithmic}
\end{algorithm}

Similar to Algorithm \ref{alg_CRAIG_G_SGSPP}, the main advantage of Algorithm \ref{alg_nsCRAIG_G_NSGSPP} is that it also does not contain references to $E$ and $F$ and works directly with the formulation (\ref{e11}).
When $C=O$, Algorithm \ref{alg_nsCRAIG_G_NSGSPP} reduces to the nsCRAIG solver for nonsymmetric saddle point systems developed in \cite{ADCKUR2025}.
Algorithm \ref{alg_nsCRAIG_G_NSGSPP} requires one more matrix-vector product $s_k=Cr_k$ and one more inner product $r_{k}^Ts_{k}$ in each iteration, and also requires two more vectors of storage: $s$ and $t$ than the nsCRAIG solver in \cite{ADCKUR2025}.
In the algorithm, it is sufficient to store only the latest left vector 
$[v_{k}; t_{k}]$ of size $m+n$ to compute $[v_{k+1}; t_{k+1}]$.
However, it is necessary to store all the right vectors $q_{k}$ of size $n$ in $Q_k$ and use them in the orthogonalization process to maintain global mutual orthogonality w.r.t. $\langle\cdot,\cdot\rangle_{N}$. The case is similar to the case of $C=O$.


Note that in step 8 of Algorithm \ref{alg_nsCRAIG_G_NSGSPP}, the vector $h_k = [h_{1,k}, h_{2,k}, \cdots, h_{k,k}]^T$ stores the inner products that are needed to orthogonalize $\hat{g}_k$ against all the previous right vectors in $Q_k$ w.r.t. $\langle\cdot,\cdot\rangle_{N}$. This is expressed by using the Gram-Schmidt method for simplicity, which is identical to that in step 8 of Algorithm \ref{alg_NS_GKB}. However, the vectors can rapidly lose orthogonality in applications where
matrices are ill-conditioned due to the accumulation of rounding error stemming from floating-point arithmetic.
In practice, the Gram-Schmidt method can be replaced with Modified Gram-Schmidt or orthogonal transformations such as Givens rotations \cite{YS2003} 
to improve numerical reliability.


Next, we give error estimates for the errors on $u-u^{(k)}$ and $p-p^{(k)}$, and on the dual norm of the residual $r^{(k)} =b - A^T u^{(k)} + Cp^{(k)}$.

Given the structure of $\beta_1e_1$ and $B_k$, let $x_k=L_k^Tz_k$ and by the first equation in (\ref{e44}), we can solve $B_k^Tx_k=\beta_1e_1$ in a recursive manner, i.e.,
\begin{equation}\label{e47}
\chi_1=\frac{\beta_1}{\alpha_1},\quad \chi_{k}=-\frac{\beta_k}{\alpha_k}\chi_{k-1},\quad L_k^Tz_k = x_k = \left[\begin{array}{c} x_{k-1}\\ \chi_{k}
\end{array}\right].
\end{equation}
The process is similar to that given in (\ref{e324}).

By (\ref{e317}), (\ref{e41}), (\ref{e44}), (\ref{e45}), (\ref{e46}) and (\ref{e47}), at step $k$ of Algorithm \ref{alg_nsCRAIG_G_NSGSPP}, we have
\begin{equation}\label{e48}
\begin{array}{l}
\|u-u^{(k)}\|_{M}^2+(p-p^{(k)})^TC(p-p^{(k)}) = \|p-p^{(k)}\|_{A^TM^{-T}A+C}^2=\|p-p^{(k)}\|_{S^T}^2
\\
\quad =\|V_nz_n-V_kz_k\|_{M}^2+(Q_ny_n-Q_ky_k)^TC(Q_ny_n-Q_ky_k)\\
 \quad=\|V_nz_n-V_kz_k\|_{M}^2+(-Q_nB_n^{-1}z_n+Q_kB_k^{-1}z_k)^TC(-Q_nB_n^{-1}z_n+Q_kB_k^{-1}z_k)\\
 \quad=\|V_nL_n^{-T}x_n-V_kL_k^{-T}x_k\|_{M}^2+(D_nL_n^{-T}x_n-D_kL_k^{-T}x_k)^TC(D_nL_n^{-T}x_n-D_kL_k^{-T}x_k)\\
 \quad=[0\ x_{n-k}^T]L_n^{-1}(V_n^TMV_n+D_n^TCD_n)L_n^{-T}\left[\begin{array}{c}0\\ x_{n-k} \end{array}\right]=[0\ x_{n-k}^T]L_n^{-T}\left[\begin{array}{c}0\\ x_{n-k} \end{array}\right]\\
 \quad=x_{n-k}^TL_{n-k}^{-T}x_{n-k} =x_{n-k}^Tz_{n-k}\\
 \quad=\sum_{i=k+1}^n\chi_i\zeta_i,\\
\end{array}
\end{equation}
where 
\begin{equation*}
L_n=\left[\begin{array}{cc}L_k & O \\ L_{n-k,k}& L_{n-k}\end{array}\right], \quad
L_n^{-T}=
\left[\begin{array}{cc}L_k^{-T} & -L_k^{-T}L_{n-k,k}^TL_{n-k}^{-T} \\ O& L_{n-k}^{-T}\end{array}\right],
\end{equation*}
with $L_k\in R^{k\times k}$, $L_{n-k,k}\in R^{(n-k)\times k}$, $L_{n-k}\in R^{(n-k)\times (n-k)}$, $z_{n}=[z_{k}; z_{n-k}]$, $x_n=[x_{k}; x_{n-k}]=L_n^Tz_n$, and $z_{n-k}=L_{n-k}^{-T}x_{n-k}=[\zeta_{k+1}, \zeta_{k+2}, \dots, \zeta_{n}]^T$  with $x_{n-k}=[\chi_{k+1}, \chi_{k+2}, \dots, \chi_{n}]^T$.
The elements in $z_{n-k}$ that is required in (\ref{e48}) can be obtained by computing $n-k$ steps in the backwards substitution with $L_{n}^{T}$.
It is easy to verify that the (\ref{e48}) is the error in the energy norm for the primal variable of (\ref{e17}), i.e., $\|\bar{u}\|_{\bar{M}}$ rather than that for the primal variable of (\ref{e11}).
Since $C$ is SPSD, by (\ref{e48}), we have 
\begin{equation*}
\|u-u^{(k)}\|_{M}^2\leq \|p-p^{(k)}\|_{A^TM^{-T}A+C}^2= \sum_{i=k+1}^n\chi_i\zeta_i.
\end{equation*}
Similar to the stopping criterion (\ref{e329}) for Algorithm \ref{alg_CRAIG_G_SGSPP}, 
the relative error in (\ref{e49})
can be used as a stopping criterion of Algorithm \ref{alg_nsCRAIG_G_NSGSPP}
\begin{equation}\label{e49}
\bar{\xi}_{k,d}^2=\frac{\sum_{i=k-d+1}^{k}\chi_i\zeta_i}{\sum_{i=1}^{k}\chi_i\zeta_i},
\end{equation}
where $x_k=[\chi_1,\dots, \chi_{k}]^T$ can be cheaply calculated by (\ref{e47}), but $z_{k}=[\zeta_1,\dots, \zeta_{k}]^T$ is obtained by solving the system with upper triangular matrix $L_{k}^{T}$ at each iteration.
The computational cost of the latter increases as $k$ increases. Next, we give a more attractive stopping criterion 
for Algorithm \ref{alg_nsCRAIG_G_NSGSPP}. 

By (\ref{e21}), (\ref{e41}),  (\ref{e44}), (\ref{e45}) and (\ref{e47}), along with the structure of $e_k$ and $L_k$, we have
\begin{equation}\label{e410}
\begin{array}{ll}
r^{(k)}_{nsCRAIG}=b-A^Tu^{(k)}+Cp^{(k)}&=\beta_1Nq_1-A^TV_kz_k+CQ_ky_k \\
&=\beta_1Nq_1-A^TV_kz_k+T_kB_ky_k \\
&= \beta_1Nq_1-(A^TV_k+T_k)z_k \\
&= \beta_1Nq_1-(NQ_kH_k+\beta_{k+1}Nq_{k+1}e_k^T)z_k \\
&= \beta_1NQ_ke_1 -\beta_1NQ_ke_1- \beta_{k+1}Nq_{k+1}e_k^Tz_k \\
&= - \beta_{k+1}Nq_{k+1}e_k^TL_k^{-T}x_k \\
&= - \beta_{k+1}Nq_{k+1}e_k^Tx_k \\
&= - \chi_k\beta_{k+1}Nq_{k+1}. \\
\end{array}
\end{equation}
It follows from (\ref{e410}) that the residual of the second equation of (\ref{e11}) is parallel to the vector $Nq_{k+1}$ and it is orthogonal to all right vectors stored in $Q_k$, i.e., $(r^{(k)})^TQ_k = 0$ owing to the $N$-orthogonality of $q_{k+1}$, like the symmetric case. Using (\ref{e410}) and the first equation in (\ref{e46}), we have the residual of (\ref{e11}) at $(u^{(k)},p^{(k)})$:
\begin{equation}\label{e411}
\mathbf{res}_k^{nsCRAIG} \triangleq
\left[\begin{array}{c}
0\\ b
\end{array}\right]-\left[\begin{array}{cc}
M& A\\
A^T & -C\\
\end{array}\right]\left[\begin{array}{c}
u^{(k)}\\ p^{(k)}
\end{array}\right]
=\left[\begin{array}{c}
O\\ r^{(k)}_{nsCRAIG}
\end{array}\right]
=\left[\begin{array}{c}
0\\ - \chi_k\beta_{k+1}Nq_{k+1}
\end{array}\right].
\end{equation}
Then,
at step $k$ of Algorithm  \ref{alg_nsCRAIG_G_NSGSPP}, we have the following residual norm
\begin{equation}\label{e412}
\|\mathbf{res}_k^{nsCRAIG}\|_{D_0^{-1}}=\|r^{(k)}_{nsCRAIG}\|_{N^{-1}}=
\|b-A^Tu^{(k)}+Cp^{(k)}\|_{N^{-1}}=\beta_{k+1}|\chi_k|,
\end{equation}
where $D_0=\mbox{blkdiag}(M,N)$ and the $N$-orthogonality of $q_{k+1}$ is used.
The residual norm (\ref{e411}) in the case of $C=O$ and $N=I_n$ was given in Lemma 3.4 in \cite{ADCKUR2025}.
With this bound for the absolute residual, we can devise one for the relative residual in (\ref{e413}), which can be used as another stopping criterion for Algorithm \ref{alg_nsCRAIG_G_NSGSPP}:
\begin{equation}\label{e413}
\frac{\|\mathbf{res}_k^{nsCRAIG}\|_{D_0^{-1}}}{\|\mathbf{res}_0^{nsCRAIG}\|_{D_0^{-1}}}
=\frac{\|r^{(k)}_{nsCRAIG}\|_{N^{-1}}}{\|r^{(0)}_{nsCRAIG}\|_{N^{-1}}}
=\frac{\|b-A^Tu^{(k)}+Cp^{(k)}\|_{N^{-1}}}{\|b\|_{N^{-1}}}=\frac{\beta_{k+1}}{\beta_{1}}|\chi_k|.
\end{equation}
Since $\beta_{k+1}$ and $\chi_k$ can be recursively calculated, the residual norm can be obtained very cheaply.
However, for challenging problems, there may be a gap between the error norm and the residual norm, making the latter less reliable as stopping criterion.
In specific numerical experiments, we still use (\ref{e413}) rather than (\ref{e49}) as the stopping rule due to its simplicity and generality, see steps 9-12 in Algorithm \ref{alg_nsCRAIG_G_NSGSPP}.

\begin{remark} \label{R41}
Following the process in Remark \ref{R31}, we also obtain that Algorithm \ref{alg_nsCRAIG_G_NSGSPP} cannot terminate with $\alpha_{k+1} = 0$ since $b\in \mbox{Range}([A^T\ C])$.
If $\beta_{k+1}=0$, then by (\ref{e411}),
Algorithm \ref{alg_nsCRAIG_G_NSGSPP} also terminates with an exact solution like Algorithm \ref{alg_CRAIG_G_SGSPP}.
\end{remark}

From (\ref{e41}) or (\ref{e42}), we can draw a link between the decomposition of the augmentation of the off-diagonal block $A$ of (\ref{e11}) or the decomposition of the off-diagonal blocks $[A^T\ C]^T$ and $[A^T\ I_n]$ of (\ref{e18}) and the decomposition of the centered-preconditioned Schur complement $N^{-\frac{1}{2}}SN^{-\frac{1}{2}}$ with $S=A^TM^{-1}A+C$, which is exactly identical to that given in (\ref{e210}).
The equation (\ref{e210}) is the matrix form of the Arnoldi process specific to FOM with zero initial guess for the NSPD systems of the form (\ref{e336}).
It is emphasized that the FOM for (\ref{e336}) does not break down since 
the coefficient matrix $N^{-\frac{1}{2}}SN^{-\frac{1}{2}}$ is positive definite.

An important property of the Arnoldi vectors in $N^{\frac{1}{2}}Q_k$ is that they lie in the Krylov subspace $\mathcal{K}_k=\mathcal{K}_k(N^{-\frac{1}{2}}SN^{-\frac{1}{2}},N^{-\frac{1}{2}}b)$ defined in (\ref{e337}).
Similar to CG, at iteration $k$ of FOM for (\ref{e336}), we look for an approximate solution $\hat{p}^{(k)}_{FOM}=N^{\frac{1}{2}}Q_ky_k$ such that
characterized variationally via
\begin{equation}\label{e414}
\hat{p}^{(k)}_{FOM}\in \mathcal{K}_k(N^{-\frac{1}{2}}SN^{-\frac{1}{2}},N^{-\frac{1}{2}}b), \quad \hat{r}^{(k)}_{FOM}=-N^{-\frac{1}{2}}b-N^{-\frac{1}{2}}SN^{-\frac{1}{2}}\hat{p}^{(k)}_{FOM}
\perp \mathcal{K}_k(N^{-\frac{1}{2}}SN^{-\frac{1}{2}},N^{-\frac{1}{2}}b).
\end{equation}
This implies that
$y_k$ satisfies that
\begin{equation}\label{e417}
H_kB_ky_k 
=-\beta_1e_1
\end{equation}
with $\beta_1=\|N^{-\frac{1}{2}}b\|_2=\|b\|_{N^{-1}}$.
Thus, by the above analysis, we have the solution of (\ref{e13}) 
at the $k$th step of FOM with preconditioner $N$:
\begin{equation*}
p^{(k)}_{FOM}=N^{-\frac{1}{2}}\hat{p}^{(k)}_{FOM}=N^{-\frac{1}{2}}N^{\frac{1}{2}}Q_ky_k=-Q_k(B_k^{-1}(H_k^{-1}\beta_1e_1)),
\end{equation*}
which is exactly the lower block $p^{(k)}$ of the approximate solution to (\ref{e11}) given in (\ref{e45}). Then, by the first equation of (\ref{e46}),  we have the following result.

\begin{theorem} 
\label{T41}
	The nsCRAIG iterates $u^{(k)}$ produced in Algorithm \ref{alg_nsCRAIG_G_NSGSPP} are related to the iterates $p^{(k)}$ of the full orthogonalization method applied to the Schur-complement equation (\ref{e13})
with preconditioner $N$ according to $u^{(k)}=-M^{-1}Ap^{(k)}$.
\end{theorem}


Due to the fact that the FOM for the NSPD systems (\ref{e13}) does not feature a minimization property, and the FOM is theoretically equivalent to our proposed nsCRAIG method according to Theorem \ref{T41}, the nsCRAIG algorithm proposed for the nonsymmetric generalized saddle point problem (\ref{e11}) also does not have a minimization property.
This implies that we cannot guarantee a decrease in the residual/error norm at each iteration of the nsCRAIG algorithm.

\begin{remark} \label{R42}
Theorem \ref{T41} establishes the important connection: for the nonsymmetric generalized saddle point problem (\ref{e11}), our proposed nsCRAIG solver and the SCR(FOM) method
produce the same iterates $p^{(k)}$ and $u^{(k)}$ until the stopping condition is fulfilled, like the case of $C=O$ given in \cite{ADCKUR2025}.
Therefore, the total computational cost of the two methods for the solution of (\ref{e11}) is exactly the same.
However, 
in the presence of ill-conditioning and if a high quality solution is required, nsCRAIG can be more accurate than FOM. This is due to the Schur complement potentially having a much higher condition number and leading to faster accumulation of errors, similar to the case regarding the CRAIG and SCR(CG) 
discussed in Remark \ref{R33}.
\end{remark}

\section{Numerical experiments} \label{Sec5}

\quad Our application  of choice is from the field of Computational Fluid Dynamics, in the form of Stokes and Navier-Stokes flow problems. To generate the linear systems for our tests, we make use of the Incompressible Flow $\&$ Iterative Solver Software\footnote{http://www.cs.umd.edu/elman/ifiss3.6/index.html} (IFISS) package, see also \cite{HCEARDJS2014,DSHEAR2014}. The particular symmetric and nonsymmetric problems under consideration are those also used in \cite{VDADCKUR2023} and \cite{HEVEHJSDSRT2008}, respectively, and described in more detail in \cite{HCEDJSAJW2014}. We briefly summarize them here.


In each case, we generate the 2D Stokes problem with IFISS, which is given by
\begin{equation}\label{e51}
\begin{array}{r}
- \Delta \vec{u}+\nabla p=0,\\
\nabla \cdot \vec{u}=0.
\end{array}
\end{equation}
We use $Q1$-$P0$ 
finite element discretization, leading to a linear system of the form (\ref{e11}) with $M \succ O$ and $C \succeq O$.
Next, we briefly describe three test problems.


%

\textbf{Test case 1: Flow over a backward facing step.} This example 
represents a flow with a parabolic inflow velocity profile passing through a domain $\Omega = ((-1, 5) \times (-1, 1)) \backslash ((-1, 0] \times (-1, 0])$. The boundary conditions are
\begin{equation}\label{e52}
\left\{
\begin{array}{ll}
u_x=4y(1-y), \quad u_y=0& \mbox{at}\  \mbox{the}\ \mbox{inflow} \ \Gamma_{in}=\{-1\}\times [0,1], \\
\frac{\partial u_x}{\partial x}-p=0,\quad \frac{\partial u_y}{\partial x}=0 & \mbox{at}\ \mbox{the}\ \mbox{outflow} \ \Gamma_{out}=\{5\}\times [-1,1], \\
\mbox{no-slip} & \mbox{on}\  \mbox{the}\ \mbox{horizontal} \ \mbox{walls}.
\end{array}
\right.
\end{equation}


After discretization, the sizes of the blocks in the resulting generalized saddle point system (\ref{e11}) are defined by $m=362498$, $n=180224$. 

\textbf{Test case 2: Driven cavity flow.} The domain for this problem is $\Omega = (-1, 1)\times(-1, 1)$, with the following boundary conditions
\begin{equation}\label{e53}
\left\{
\begin{array}{ll}
u_x=1-x^4,\quad u_y=0& \mbox{on}\  \mbox{the}\ \mbox{wall} \ \Gamma_{top}=[-1,1]\times \{1\}, \\
\mbox{no-slip} & \mbox{on}\  \mbox{the}\ \mbox{bottom}\ \mbox{and}\ \mbox{vertical} \ \mbox{walls}.
\end{array}
\right.
\end{equation}
This represents a model where the cavity lid is moving according to the given regularized condition, driving the enclosed flow.

After discretization, the (1,2)-block of the generalized saddle point matrix generated is rank deficient.
Thus, its first two columns are dropped and the first two rows and columns of the (2,2)-block of the matrix are dropped correspondingly, then the resulting system is the generalized saddle point systems (\ref{e11}). The sizes of the blocks in (\ref{e11}) are defined by
$m=132098$, $n=65536-2=65534$.


\textbf{Test case 3: Poiseuille flow problem.} This problem is a steady Stokes problem with the exact solution 
\begin{equation}\label{e54}
u_x=1-y^2, \quad u_y=0, \quad p=-2x+\mbox{constant}.
\end{equation}
The boundary conditions are the same as those given in 
(\ref{e52}).

After discretization, the sizes of the blocks in the resulting generalized saddle point system (\ref{e11}) are defined by $m = 105666$ and $n = 51200$. The length of the channel domain for this test case is chosen as a large number $L=1024$. Thus, the off-diagonal block $A$ of (\ref{e11}) is very likely ill-conditioned \cite{ADCKUR2024} and then the related Schur complement $S=A^TM^{-1}A+C$ can be more ill-conditioned
since the condition number of $S$ is higher than that of $A$. See \cite{EVCMAO2000} for a related analysis on that the ratio between the length and width of the channel impacts the condition number of the Schur complement.

As nonsymmetric test problem, we will use the nonlinear Navier-Stokes problem generated by IFISS in each case, which is given by
\begin{equation}\label{e55}
\begin{array}{r}
-\nu \nabla^2 \vec{u}+\vec{u}\cdot\nabla \vec{u} +\nabla p=\vec{f},\\
\nabla \cdot \vec{u}=0,
\end{array}
\end{equation}
with the kinematic viscosity $\nu$. To deal with the nonlinearity of the Navier-Stokes equation in the convection term, Picard's iteration is used to obtain the following linearized equations
\begin{equation}\label{e56}
\begin{array}{rr}
-\nu \Delta \vec{u}^{(k)}+(\vec{u}^{(k-1)}\cdot\nabla) \vec{u}^{(k)} +\nabla p^{(k)}=\vec{f}, & \mbox{in} \ \Omega, \\
\nabla \cdot \vec{u}^{(k)}=0,& \mbox{in} \ \Omega,
\end{array}
\end{equation}
for each iteration $k$, starting from an arbitrary initial guess $(\vec{u}^{(0)}, p^{(0)})$. For the linearized system (\ref{e56}), we consider a $Q1$-$P0$ 
finite element discretization, leading to a linear system of the form (\ref{e11}) 
with $M$ NSPD and $C\succeq O$. 
Next, we briefly describe three test problems.

%
%

\textbf{Test case 4: Flow over a backward facing step.} This case represents a flow with a parabolic inflow velocity profile passing through a domain $\Omega = ((-1, 5) \times (-1, 1)) \backslash ((-1, 0] \times (-1, 0])$. The boundary conditions at the inflow $\Gamma_{in}$ and on the horizontal walls are exactly the same as that given in (\ref{e52}). Only the boundary condition at the outflow $\Gamma_{out}$ changes into
\begin{equation}\label{e57}
\nu \frac{\partial u_x}{\partial x}-p=0,\quad \frac{\partial u_y}{\partial x}=0.
\end{equation}

After discretization, the sizes of the blocks in the resulting generalized saddle point system (\ref{e11}) are defined by
$m=91138$, $n=45056$. The viscosity parameter for this test case is $\nu = 1/1000$.

\textbf{Test case 5: Driven cavity flow.} The domain for this problem is $\Omega = (-1, 1)\times(-1, 1)$, with the boundary conditions given as in (\ref{e53}).

After discretization, the (1,2)-block of the block $2\times 2$ system generated is rank deficient. By dropping the first two columns of the (1,2)-block and the first two rows and columns of the (2,2)-block of the system, the resulting system is the generalized saddle point system (\ref{e11}). The block sizes in (\ref{e11}) are 
$m=33282$, $n=16384-2=16382$. The viscosity parameter for this test case is $\nu = 1/1000$.



\textbf{Test case 6: Poiseuille flow problem.} This problem is a Navier-Stokes equation with the exact solution
\begin{equation}\label{e58}
u_x=1-y^2,\quad u_y=0,\quad p=-2\nu x+\mbox{constant}.
\end{equation}
The only difference between (\ref{e58}) and (\ref{e54}) is that
the pressure gradient is proportional to the viscosity parameter.
The boundary conditions at the inflow $\Gamma_{in}$ and on the horizontal walls are given in (\ref{e52}), and the condition at the outflow $\Gamma_{out}$ is given in (\ref{e57}).


After discretization, the block sizes of the resulting generalized saddle point system (\ref{e11}) are defined by $m=27234$, $n=12800$. The length of the channel domain and viscosity parameter for this test case are $L=1024$ and $\nu = 1/1000$, respectively. In this case, since the chosen $L$ is large, the off-diagonal block of (\ref{e11}) and the associated Schur complement are ill-conditioned similar to the symmetric case.

For the symmetric generalized saddle point problems (\ref{e11}), the CRAIG algorithm proposed in Section \ref{Sec3} is theoretically equivalent the SCR
method with inner CG iteration on the associated Schur-complement equation with an SPD preconditioner $N$.
Similarly, in Section \ref{Sec4}, for the nonsymmetric generalized saddle point problems (\ref{e11}), the proposed nsCRAIG algorithm is theoretically equivalent to the SCR
method in which FOM iteration is applied to the preconditioned Schur-complement equation of the problem.
We firstly test whether these equivalences also hold numerically.

Generalized saddle point systems are often solved with the MINRES method \cite{CCPMAS1975} when the leading block is symmetric and with the GMRES method \cite{YSMHS1986} when this condition is not fulfilled.
These methods belong to the class of coupled algorithms compared to our proposed segregated algorithms.
We secondly test whether the performance of our proposed segregated CRAIG and nsCRAIG solvers outperform that of the common coupled MINRE and GMRES methods.

It is crucial that all solvers have comparable costs as we compared the number of iterations performed by each solver. In this regard, the most expensive operations in CG, FOM, and our proposed CRAIG and nsCRAIG algorithms are to apply $M^{-1}$ and $N^{-1}$ to vectors.
To bring MINRES to the same level as CRAIG, we consider it as applied to the centered-preconditioned problem
\begin{equation}\label{e59}
\left[\begin{array}{c c}
		M & \\
	      & N\\
	\end{array}\right]^{-\frac{1}{2}}\left[\begin{array}{c c}
		M & A\\
		A^T & -C\\
	\end{array}\right]\left[\begin{array}{c c}
		M & \\
	      & N\\
	\end{array}\right]^{-\frac{1}{2}}
	\left[\begin{array}{c}
		\tilde{u} \\ \tilde{p}
	\end{array}\right]
	=\left[\begin{array}{c c}
		M & \\
	      & N\\
	\end{array}\right]^{-\frac{1}{2}}\left[\begin{array}{c}
		0 \\ b
	\end{array}\right],
\end{equation}
\begin{equation*}
\left[\begin{array}{c}
		u \\ p
	\end{array}\right]=\left[\begin{array}{c c}
		M & \\
	      & N\\
	\end{array}\right]^{-\frac{1}{2}}
	\left[\begin{array}{c}
		\tilde{u} \\ \tilde{p}
	\end{array}\right].
\end{equation*}
Similarly, to bring GMRES to the same level as nsCRAIG, we consider it as applied to the right-preconditioned problem
\begin{equation}\label{e510}
\left[\begin{array}{c c}
		M & A\\
		A^T & -C\\
	\end{array}\right]\left[\begin{array}{c c}
		M & \\
	      & N\\
	\end{array}\right]^{-1}
	\left[\begin{array}{c}
		\tilde{x} \\ \tilde{y}
	\end{array}\right]
	=\left[\begin{array}{c}
		0 \\ b
	\end{array}\right],\quad \left[\begin{array}{c}
		x \\ y
	\end{array}\right]=\left[\begin{array}{c c}
		M & \\
	      & N\\
	\end{array}\right]^{-1}
	\left[\begin{array}{c}
		\tilde{x} \\ \tilde{y}
	\end{array}\right].
\end{equation}
Since GMRES treats the matrix system as a whole, i.e., in an all-at-once manner, it
stores vectors of length $m+n$ in memory. In contrast, only the right vectors of size $n$ need to be stored in our proposed nsCRAIG. For the specific problems in the test cases 4-6 and the $Q1$-$P0$ finite element discretization, $n$ is approximately half of $m$.
Hence, for a given amount of memory, the proposed nsCRAIG solver can perform more iterations than GMRES.

For all the experiments, we report the summary of the results obtained using a MATLAB version of algorithms, where the matrix $M$ is factorized using the MATLAB function \texttt{chol} $([R, flag, p] = chol(M,`vector'))$ if $M$ is SPD and \texttt{lu} $([L, U, p, q] = lu(M,`vector'))$ if $M$ is NSPD.
Thus, every inner problem of applying $M^{-1}$ to a vector is solved exactly, which can stay close to the theory and enable our solvers to continue converging.
Efforts to solve this inner problem of the generalized GKB in an approximate fashion with an iterative solver have only been explored for the case of $M \succ O$ and $C=O$ in \cite{VDADCKUR2023}.

Apparently, the convergence of all the considered
solvers depend on the choice of the SPD preconditioner $N$.
Using the equivalence given in Theorems \ref{T31} and \ref{T41} and the explanation in \cite[Section 10.1.1]{MBGHGJL2005},
$N$ should be chosen as a good preconditioner for the Schur complement $S=A^TM^{-1}A+C$ of the system (\ref{e11}).
As explained in \cite{HCEDJSAJW2014}, for Stokes problems, a good choice for $N$ is the pressure mass matrix $Q$. The choice $N=(1/\nu)Q$ also has merit for discrete Navier-Stokes problems for low Reynolds number. Here, the Reynolds number is a quantity inversely proportional to $\nu$. However, it does not take into account the effects of convection on the Schur complement operator, and convergence rates deteriorate as the Reynolds number increases. The pressure convection-diffusion and least-squares commutator preconditioners of the Schur complement $S$ is ideal choices for $N$ that better reflect the balance of convection and diffusion in the problem and so lead to improved convergence rates at higher Reynolds numbers. In this section, for simplicity, we choose $N=Q$ and $N=(1/\nu)Q$ for discrete systems from Stokes and Navier-Stokes problems, respectively.
Note that in the case of $Q1$-$P0$ finite element discretization, the matrix $Q$ is diagonal, so
the application of $N^{-1}$ to a vector is cheaply by multiplying scalar quantities.

In the section, we compare the proposed CRAIG for (\ref{e11}) 
to MINRES for (\ref{e57}) and compare the nsCRAIG for (\ref{e11}) to
GMRES for (\ref{e58}) by measuring the necessary number of iterations, the elapsed CPU times in seconds, and 
the relative error norm (denoted by ``ERR") to reach convergence, defined as reducing the relative residual norm (denoted by ``RES") below the given tolerance $tol$ or exceeding the specified maximum iteration number $k_{\max}=3000$.
Here,
\begin{equation*}
\mbox{ERR} \triangleq \frac{\|z^{(k)}-z^*\|_2}{\|z^{(0)}-z^*\|_2}, \quad \mbox{RES} \triangleq \frac{\|f-\mathcal{A}z^{(k)}\|_2}{\|f-\mathcal{A}z^{(0)}\|_2},
\end{equation*}
where $z^{(k)}=[u^{(k)};p^{(k)}]$ is the $k$th approximate solution of the tested linear systems (\ref{e11}). Note that the ``RES" of the proposed CRAIG and nsCRAIG solvers are identical to (\ref{e333}) and (\ref{e413}), respectively.
Besides, the corresponding right-hand-side vectors 
such that the exact solutions of the tested problems are $z^*=[1,1,\cdots,1]^T \in \mathbb{R}^{m+n}$. All the initial vectors for the MINRES for (\ref{e59}) and GMRES for (\ref{e510}) are set to be zero, and all the initial vectors for the inner CG and FOM iterations in the SCR methods are also set to be zero.
The tolerances $tol$ chosen for all the test cases are $10^{-6}$ and $10^{-15}$.
All of the experiments are performed in this paper using MATLAB (version 24.1.0.2578822 (R2024a)) on a PC equipped with 13th Gen Intel(R) Core(TM) i5-13600KF 3.50 GHz, 32.0 GB RAM, and Win11 operating system.

\begin{figure}[!ht]
	\centering
\subfigure[$tol=10^{-6}$]
{\includegraphics[width=0.49\textwidth, trim=0cm 0cm 1cm 3cm, clip]{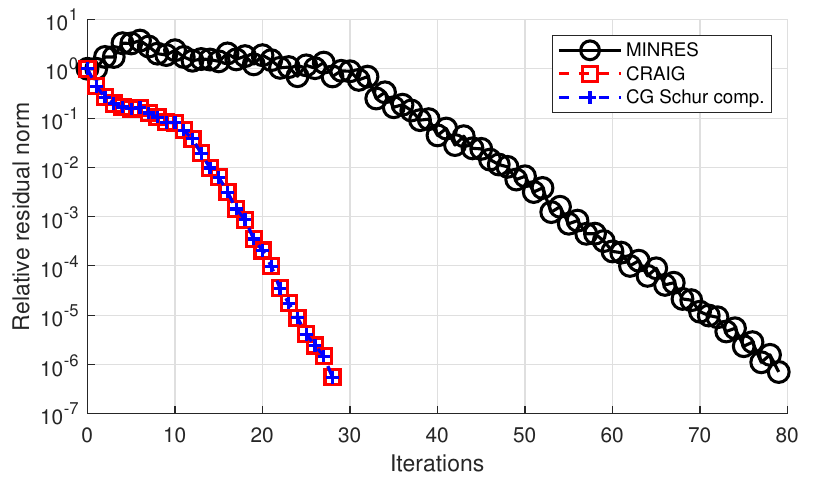}} \hfill  \subfigure[$tol=10^{-15}$]
{\includegraphics[width=0.49\textwidth, trim=0cm 0cm 1cm 3cm, clip]{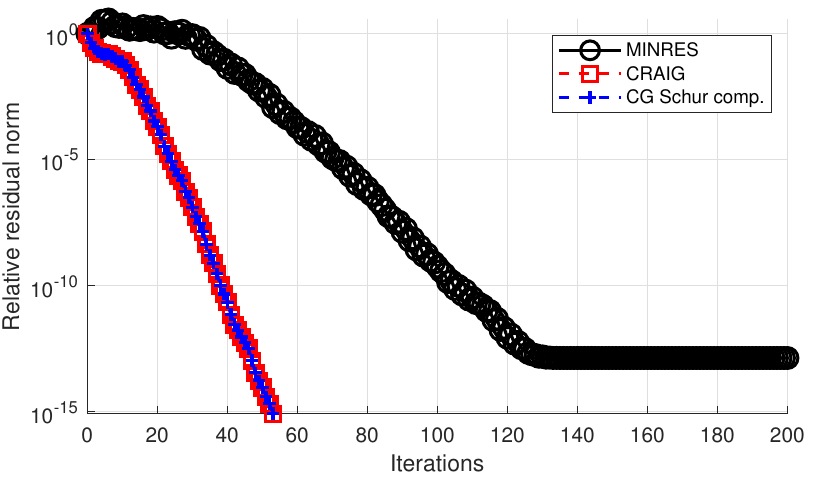}}\\
	\caption{Convergence history of the listed solvers for the Stokes flow over a backward facing step. 
}\label{f51}
\end{figure}

\begin{table}[!ht]
\caption{Numerical results of the linear solvers MINRES, CRAIG and SCR(CG)  in the case of Stokes flow over a backward facing step with $tol=10^{-6}$ and $tol=10^{-15}$. 
}\label{t51}
\begin{center}
\begin{tabular}{c|ccc|ccc}
\hline
Methods &  & $tol=10^{-6}$ &  &  & $tol=10^{-15}$  & \\
\cline{2-7}
 & MINRES& CRAIG& SCR(CG)& MINRES& CRAIG& SCR(CG) \\
\hline
iterations &79 &28 &28 &- &53 &53 \\
time &6.6172 &\textbf{1.9823} &2.0003 &- &\textbf{3.5830} &3.6524 \\
ERR &2.9820e-09 &1.3827e-07 &1.3827e-07 &- &4.9175e-12 &4.9173e-12 \\
\hline
\end{tabular}
\end{center}
\end{table}

\begin{figure}[!ht]
	\centering
\subfigure[$tol=10^{-6}$]
{\includegraphics[width=0.49\textwidth, trim=0cm 0cm 1cm 3cm, clip]{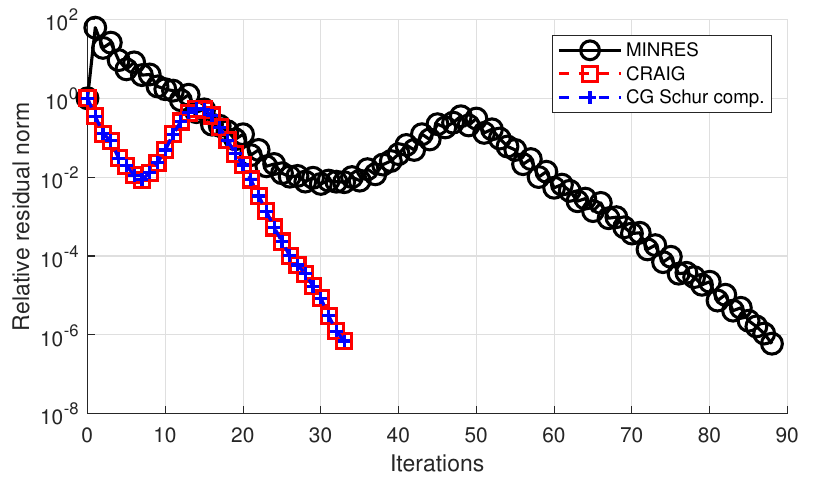}} \hfill  \subfigure[$tol=10^{-15}$]
{\includegraphics[width=0.49\textwidth, trim=0cm 0cm 1cm 3cm, clip]{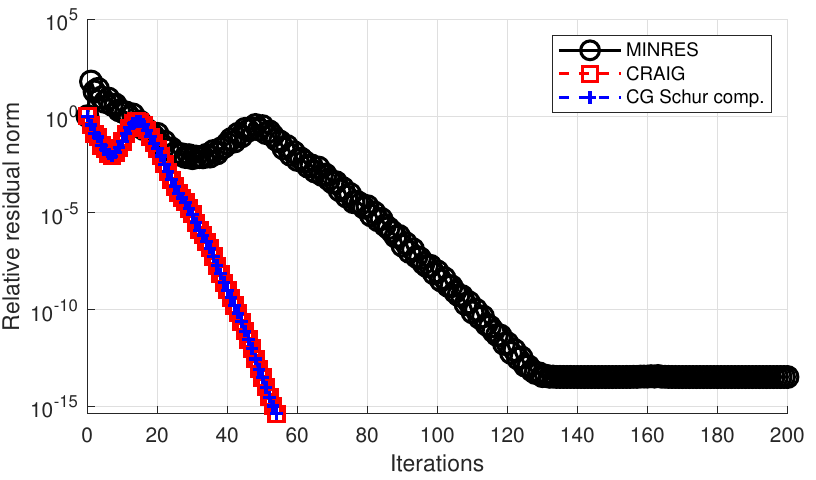}}\\
	\caption{Convergence history of the listed solvers for the Stokes driven cavity flow
.}\label{f52}
\end{figure}

\begin{table}[!ht]
\caption{Numerical results of the linear solvers MINRES, CRAIG and SCR(CG) in the case of Stokes driven cavity flow with $tol=10^{-6}$ and $tol=10^{-15}$. 
}\label{t52}
\begin{center}
\begin{tabular}{c|ccc|ccc}
\hline
Methods &  & $tol=10^{-6}$ &  &  & $tol=10^{-15}$  & \\
\cline{2-7}
 & MINRES& CRAIG& SCR(CG)& MINRES& CRAIG& SCR(CG) \\
\hline
iterations &88 &33 &33 &- &54 &54 \\
time &2.7248 &\textbf{0.8427} &0.8911 &- &\textbf{1.3497} &1.4224 \\
ERR &5.5118e-11 &1.8637e-09 &1.8637e-09 &- &5.3560e-11 &5.3557e-11 \\
\hline
\end{tabular}
\end{center}
\end{table}

\begin{figure}[!ht]
	\centering
\subfigure[$tol=10^{-6}$]
{\includegraphics[width=0.49\textwidth, trim=0cm 0cm 1cm 3cm, clip]{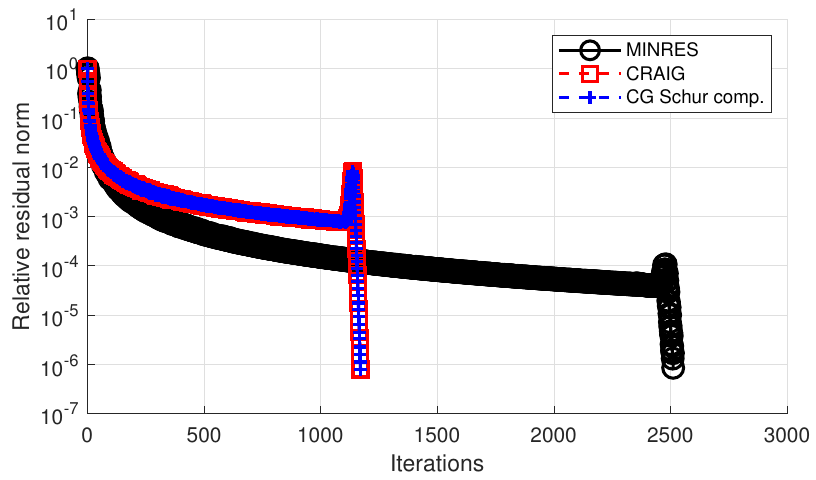}} \hfill  \subfigure[$tol=10^{-15}$]
{\includegraphics[width=0.49\textwidth, trim=0cm 0cm 1cm 3cm, clip]{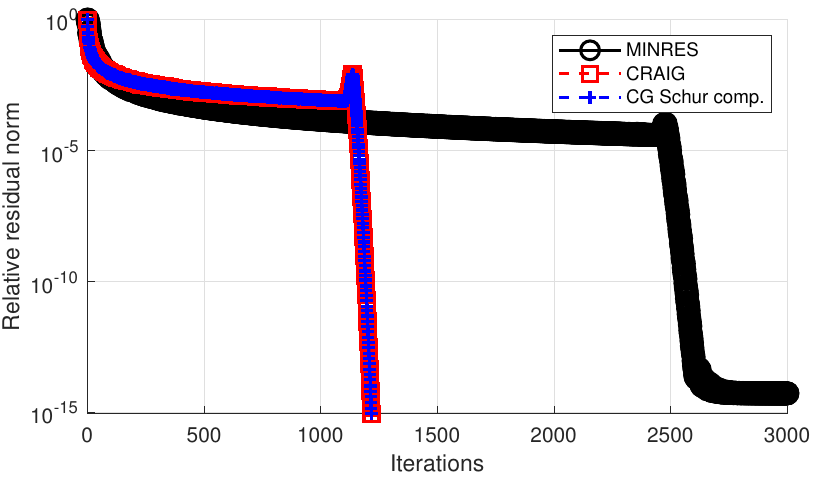}}\\
	\caption{Convergence history of the listed solvers for the Stokes flow over a channel domain
.}\label{f53}
\end{figure}

\begin{table}[!ht]
\caption{Numerical results of the linear solvers MINRES, CRAIG and SCR(CG) in the case of Stokes flow over a channel domain with $tol=10^{-6}$ and $tol=10^{-15}$. 
}\label{t53}
\begin{center}
\begin{tabular}{c|ccc|ccc}
\hline
Methods &  & $tol=10^{-6}$ &  &  & $tol=10^{-15}$  & \\
\cline{2-7}
 & MINRES& CRAIG& SCR(CG)& MINRES& CRAIG& SCR(CG) \\
\hline
iterations &2510 &1170 &1170 &- &1217 &1217 \\
time &39.7865 &14.3022 &\textbf{13.8578} &- &14.9195 &\textbf{14.3662} \\
ERR &1.1182e-05 &3.5765e-08 &3.5618e-08 &- &2.5675e-12 &2.5838e-12 \\
\hline
\end{tabular}
\end{center}
\end{table}

For the flow over a backward-facing step, the driven cavity flow and the flow over a channel domain,
we plot the convergence history of the linear symmetric solvers
in Figure \ref{f51}, Figure \ref{f52}, and Figure \ref{f53}, respectively, and tabulate their iteration counts, timings, and relative error norm in Table \ref{t51}, Table \ref{t52}, and Table \ref{t53}, respectively. The symbol ``-" in all the tables indicates that the number of iterations of the solver exceeds the given maximum iteration number $k_{max}$.
It is visible how the SCR(CG) and our CRAIG behave identically and significantly faster than the centered-preconditioned MINRES regardless of whether a moderately approximate solution or a high quality solution is required.
In terms of iterations, MINRES needs about 2-3 times more iterations
than our CRAIG algorithm to reach convergence if a moderately approximate solution is required, and MINRES needs about 2 times more iterations than our algorithm to reach convergence for test case 1 if a high quality solution is required. However, MINRES 
cannot meet the termination criterion for test cases 2 and 3 if a high quality solution is required. This may be because the resulting generalized saddle point problems from the test cases 2 and 3 are very ill-conditioned.
In terms of times, CRAIG is slightly faster than CG,
about 3-4 times faster than MINRES if a moderately approximate solution is required, while CRAIG is slightly slower than CG if a high quality solution is required. The reason is the number of iteration is relatively large, see also Remark \ref{R33}.
In terms of relative error norm, MINRES is more accurate than CRAIG for test cases 1 and 2, while it is less accurate than CRAIG for test case 3 if a moderately approximate solution is required.

It is interesting to note that the centered-preconditioned MINRES can exhibit some form of stagnation at every other step, i.e., after one step where the global residual of the system decreases, the successive step does not do it significantly.
In fact, the similar behavior for the symmetric saddle point problems 
has been noted in \cite{MA2013}.
The explanation is related to the particular choice of block preconditioner used in (\ref{e59}) that leads to a matrix with a symmetric spectrum. According to the results in \cite{BFARDJSAJW1998}, such a spectrum has an impact on the convergence of the solvers, which reduces the residual only every other step.


\begin{figure}[!ht]
	\centering
\subfigure[$tol=10^{-6}$]
{\includegraphics[width=0.49\textwidth, trim=0cm 0cm 1cm 3cm, clip]{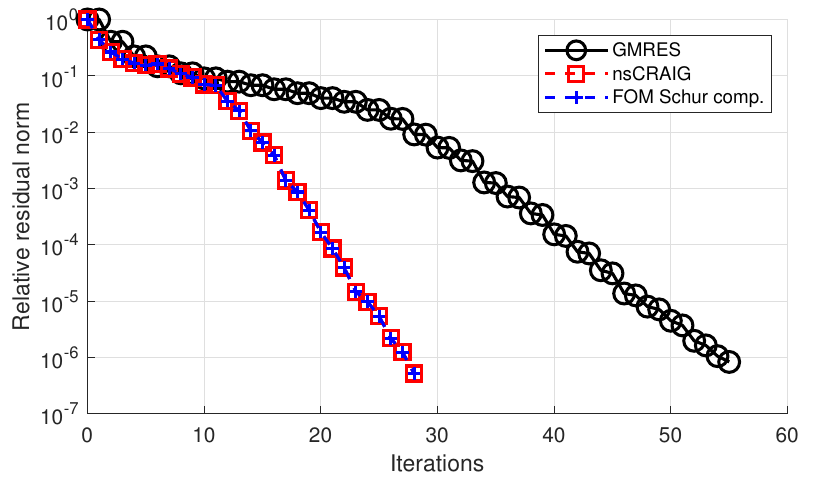}} \hfill  \subfigure[$tol=10^{-15}$]
{\includegraphics[width=0.49\textwidth, trim=0cm 0cm 1cm 3cm, clip]{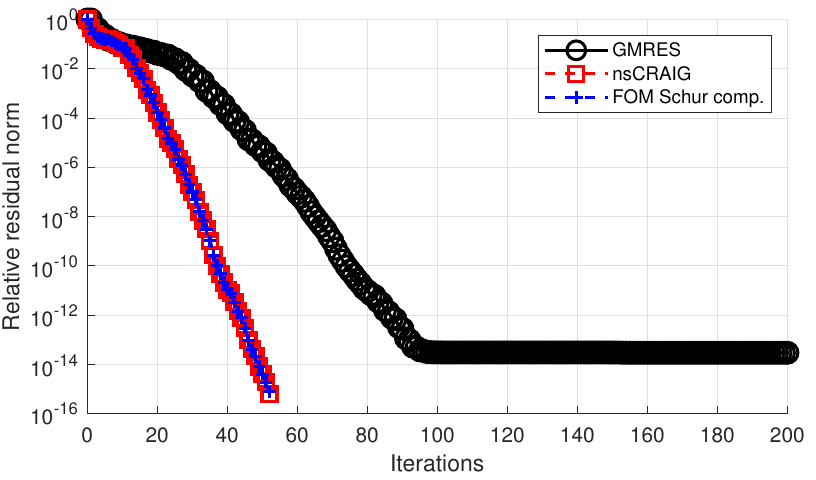}}\\
	\caption{Convergence history of the listed solvers for the Navier-Stokes flow over a backward facing step
.}\label{f54}
\end{figure}

\begin{table}[!ht]
\caption{Numerical results of the linear solvers GMRES, nsCRAIG and SCR(FOM) in the case of Navier-Stokes flow over a backward facing step  with $tol=10^{-6}$ and $tol=10^{-15}$. 
}\label{t54}
\begin{center}
\begin{tabular}{c|ccc|ccc}
\hline
Methods &  & $tol=10^{-6}$ &  &  & $tol=10^{-15}$  & \\
\cline{2-7}
 & GMRES& nsCRAIG& SCR(FOM)& GMRES& nsCRAIG& SCR(FOM) \\
\hline
iterations &55 &28 &28 &- &52 &52 \\
time &1.7760 &\textbf{0.4822} &0.5404 &- &\textbf{0.9567} &1.1668 \\
ERR &2.6140e-07 &9.3276e-08 &9.3276e-08 &- &6.8294e-13 &6.8180e-13 \\
\hline
\end{tabular}
\end{center}
\end{table}

\begin{figure}[!ht]
	\centering
\subfigure[$tol=10^{-6}$]
{\includegraphics[width=0.49\textwidth, trim=0cm 0cm 1cm 3cm, clip]{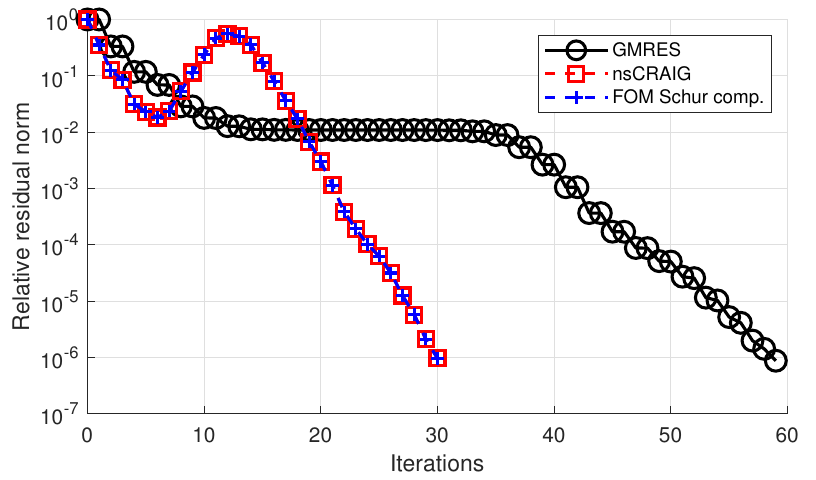}} \hfill  \subfigure[$tol=10^{-15}$]
{\includegraphics[width=0.49\textwidth, trim=0cm 0cm 1cm 3cm, clip]{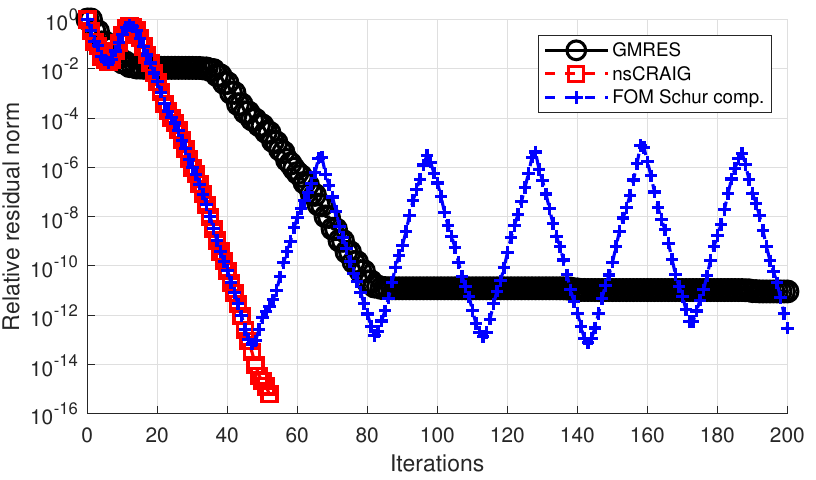}}\\
	\caption{Convergence history of the listed solvers for the Navier-Stokes driven cavity flow
.}\label{f55}
\end{figure}

\begin{table}[!ht]
\caption{Numerical results of the linear solvers GMRES, nsCRAIG and SCR(FOM) in the case of Navier-Stokes driven cavity flow with $tol=10^{-6}$ and $tol=10^{-15}$. 
}\label{t55}
\begin{center}
\begin{tabular}{c|ccc|ccc}
\hline
Methods &  & $tol=10^{-6}$ &  &  & $tol=10^{-15}$  & \\
\cline{2-7}
 & GMRES& nsCRAIG& SCR(FOM)& GMRES& nsCRAIG& SCR(FOM) \\
\hline
iterations &59 &30 &30 &- &52 &- \\
time &0.4811 &\textbf{0.1882} &0.2022 &- &0.3588 &- \\
ERR &7.2029e-09 &5.2778e-09 &5.2778e-09 &- &7.5450e-13 &- \\
\hline
\end{tabular}
\end{center}
\end{table}

\begin{figure}[!ht]
	\centering
\subfigure[$tol=10^{-6}$]
{\includegraphics[width=0.49\textwidth, trim=0cm 0cm 1cm 3cm, clip]{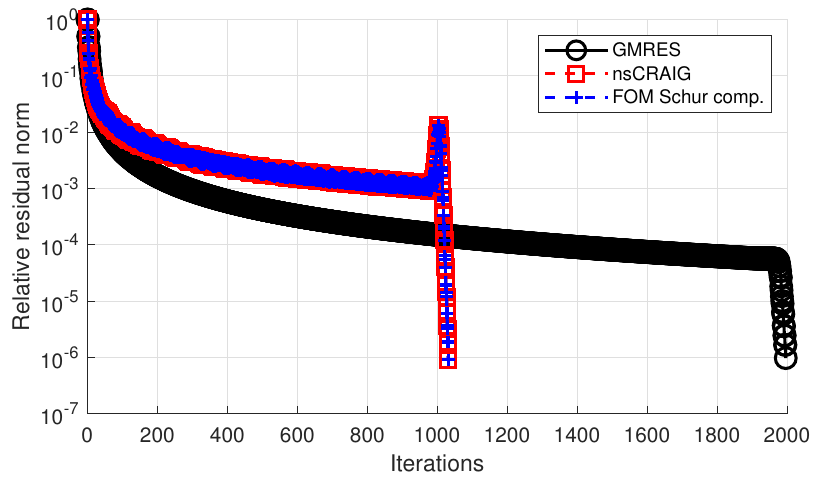}} \hfill  \subfigure[$tol=10^{-15}$]
{\includegraphics[width=0.49\textwidth, trim=0cm 0cm 1cm 3cm, clip]{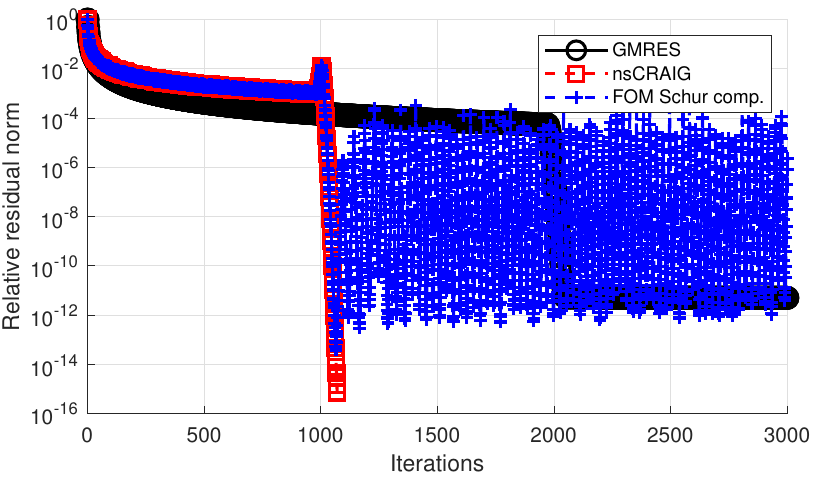}}\\
	\caption{Convergence history of the listed solvers for the Navier-Stokes flow over a channel domain
.}\label{f56}
\end{figure}

\begin{table}[!ht]
\caption{Numerical results of the linear solvers GMRES, nsCRAIG and SCR(FOM) in the case of Navier-Stokes flow over a channel domain with $tol=10^{-6}$ and $tol=10^{-15}$. 
}\label{t56}
\begin{center}
\begin{tabular}{c|ccc|ccc}
\hline
Methods &  & $tol=10^{-6}$ &  &  & $tol=10^{-15}$  & \\
\cline{2-7}
 & GMRES& nsCRAIG& SCR(FOM)& GMRES& nsCRAIG& SCR(FOM) \\
\hline
iterations & 1995 & 1031 & 1031 & - & 1070 & - \\
time & 124.1117 & \textbf{27.2722} & 40.1988 & - & 28.4551 & - \\
ERR & 1.1491e-04 & 3.3546e-08 & 3.3546e-08 & - & 5.7348e-13 & - \\
\hline
\end{tabular}
\end{center}
\end{table}


For the linearized Navier-Stokes flow over a backward-facing step, a driven cavity, and a channel domain, we plot the convergence history of
the linear nonsymmetric solvers in Figures \ref{f54}, \ref{f55}, and \ref{f56}, respectively, and tabulate the iteration counts, timings, and the relative error norms in Tables \ref{t54}, \ref{t55}, and \ref{t56}, respectively. In order to enable consistent results and avoid loss of orthogonality, all the nonsymmetric solvers we compare make use of the modified Gram-Schmidt algorithm to generate sets of mutually orthogonal vectors.
It is visible how the SCR(FOM) and our nsCRAIG behave identically when a moderately approximation solution of (\ref{e11}) is required. The SCR(FOM) and our nsCRAIG behave identically for the test case 4 when a high quality solution of (\ref{e11}) is required. However, for the test cases 5 and 6 (i.e., in the presence of ill-conditioning) and if a high quality solution is required, they behave identically except in the near convergence stage, and nsCRAIG converges successfully but SCR(FOM) fails to meet the accuracy requirements as the convergence curve of inner FOM iteration oscillates near convergence.
Whether it is a solution that requires moderate accuracy or a solution that requires high accuracy, our nsCRAIG is significantly faster than right-preconditioned GMRES.

In terms of iterations, GMRES needs about 2 times more iterations than our nsCRAIG algorithm to reach convergence with a looser stopping tolerance.
However, GMRES can not reach convergence if a tighter stopping tolerance is used.
In terms of times, if a moderately accurate solution
is required, nsCRAIG is faster than FOM, and it is anout 3-4 times faster than GMRES.
In terms of relative error norm, GMRES is less accurate than nsCRAIG and SCR(FOM) whether a moderately accurate solution or a high accurate solution is required, and GMRES and SCR(FOM) are less accurate than nsCRAIG for test cases 5 and 6 if a high accurate solution is required.

Similar to the symmetric case, the preconditioned GMRES also has phases where only every other iteration significantly contributes to the objective of reducing the residual norm. The behavior also has been noted for the nonsymmetric saddle point problems in \cite{ADCKUR2025}. The explanation is that the particular choice of block preconditioner used (\ref{e510}) makes the spectrum have its complex values distributed in a symmetric manner on both sides of a vertical line, which has an impact on the convergence of the solver \cite{BFARDJSAJW1998}.

The above experimental results indicate that in a numerical setting,
CRAIG and SCR(CG) are equivalent. Moreover, if a moderately accurate approximation is required, nsCRAIG and SCR(FOM) are also equivalent.
However, in the presence of ill-conditioning and if a high quality solution is required, nsCRAIG is more accurate 
than SCR(FOM).
This is due to the Schur complement potentially having a much higher condition number and leading to faster accumulation of errors.
Moreover, the above experimental results also show that no matter a moderately accurate or high accurate approximation is required, the proposed segregated CRAIG and nsCRAIG solvers outperform the common coupled MINRES and GMRES methods in terms of iterations and timings, and our nsCRAIG solver is more numerically stable than GMRES.

\section{Conclusions}\label{Sec6}

\quad In this paper, we have separately extended the existing generalized CRAIG solver \cite{MA2013} based on GKB for symmetric saddle point systems and nsCRAIG solver \cite{ADCKUR2025} for nonsymmetric saddle point systems to the symmetric and nonsymmetric generalized saddle point problems (\ref{e11}). 
From the theoretical point of view, the proposed CRAIG is equivalent to SCR(CG) by the known equivalence for the case of $M$ SPD and $C=O$, between generalized CRAIG \cite{MA2013} and CG. Similarly, the proposed nsCRAIG is equivalent to SCR(FOM) by the known equivalence for the case of $M$ NSPD and $C=O$, between nsCRAIG \cite{ADCKUR2025} and FOM.
Aside from the theoretical point of view, we have also illustrated this relationship in a numerical setting by our experiments.
Whether a moderately accurate or high accurate approximation is required, CRAIG and SCR(CG) are equivalent numerically.
If a moderately accurate approximation is required, nsCRAIG and SCR(FOM) are also equivalent numerically.
However, in the presence of ill-conditioning and if a high quality solution is required, nsCRAIG and SCR(FOM) are not equivalent at the stage where the algorithms are approaching convergence.

Along with our solvers' description and algorithm steps, we also provided their stopping criteria similar to that studied in \cite{MA2013,ADCKUR2025}. One is an estimate of the error for the solution of the system (\ref{e11}) in an energy norm. The other choice is an inexpensive way to compute the residual norm for (\ref{e11}) that is identical to the residual norm for the second equation of (\ref{e11}).

As we all know that MINRES and GMRES are popular choices for tackling symmetric and nonsymmetric indefinite problems, respectively. Consequently, by experiments, we compare the proposed CRAIG with MINRES and also compare the nsCRAIG with GMRES. We found that CRAIG is at least twice as fast compared to MINRES with a block diagonal preconditioner and the same is true for nsCRAIG compared to GMRES with a block diagonal preconditioner. This is due to the convergence behavior of the block-diagonal preconditioned MINRES and GMRES, where often only every other iteration significantly progresses towards convergence \cite{MA2013,ADCKUR2025}.
Moreover, similar to nsCRAIG proposed in \cite{ADCKUR2025},
our proposed nsCRAIG solver generates and stores a right basis with shorter vectors compared to GMRES and do not need to store the left basis (with long vectors), which significantly decreases memory costs.

A possible advantage of our solvers has over the equivalent SCR(CG) and  SCR(FOM) method 
is that they deliver a second basis, which corresponds to the space associated with the primal solution variable. The two bases can be explored to identify spectral information, which is useful when applying deflation. The latter mechanism can accelerate convergence for problems where the spectral distribution features outliers. Such strategies have been studied in \cite{ADCKUR2024} for the generalized CRAIG \cite{MA2013}. Similar developments in symmetric and nonsymmetric generalized saddle point problems are considered for future research.

In this paper, we have only considered exact matrix vector products of type $M^{-1}v$.
A more general and practically motivated alternative is to consider an inexact approach by making use of an iterative solver for this inner problem. For the generalized CRAIG \cite{MA2013}, such a strategy has been explored in \cite{VDADCKUR2023} and yielded promising results. A similar study concerning the methods we presented in this paper could constitute an interesting direction for future developments.



\section*{Competing interests}

The authors declare no competing interests.

\section*{Acknowledgements}

The work was supported by the Fundamental Research Funds for the Central Universities, CHD (300102123105), and National Natural Science Foundation of Shaanxi (2025JC-YBQN-086), and by National Natural Science Foundation of China (12171384).


\begin{thebibliography}{10}

\bibitem{HSDNIMGWHASAJW2006}
H. S. Dollar, N. I. M. Gould, W. H. A. Schilders, and A. J. Wathen, Implicit-factorization preconditioning and iterative solvers for regularized saddle-point systems, SIAM J. Matrix Anal. Appl., 28 (2006), pp. 170--189, https://doi.org/10.1137/05063427X.


\bibitem{DDSDO2021}
D. D. Serafino and D. Orban, Constraint-preconditioned Krylov solvers for regularized saddle-point systems, SIAM J. Sci. Comput., 43 (2021), pp. A1001--A1026, https://doi.org/10.1137/19M1291753.

\bibitem{MA2013}
M. Arioli, Generalized Golub-Kahan bidiagonalization and stopping criteria, SIAM J Matrix Anal. Appl., 34 (2013),  pp. 571--92, https://doi.org/10.1137/120866543.



\bibitem{VDADCKUR2023}
V. Darrigrand, A. Dumitrasc, C. Kruse, and U. R\"{u}de, Inexact inner-outer Golub-Kahan bidiagonalization method: A relaxation strategy, Numer. Linear Algebra Appl., 30 (2023), e2484, https://doi.org/10.1002/nla.2484.

\bibitem{DOMA2017}
D. Orban, and M. Arioli, Iterative solution of symmetric quasi-definite linear systems, SIAM, 2017, https://doi.org/10.1137/1.9781611974737.


\bibitem{ADCKUR2025}
A. Dumitrasc, C. Kruse, and U. R\"{u}de, Generalized Golub-Kahan Bidiagonalization for Nonsymmetric Saddle-Point Systems, SIAM J. Matrix Anal. Appl., 46 (2025), pp. 370--392, https://doi.org/10.1137/23M160760X.



\bibitem{HEVEHJSDSRT2008}
H. Elman, V. E. Howle, J. Shadid, D. Silvester, and R. Tuminaro, Least squares preconditioners for stabilized discretizations of the Navier-Stokes equations, SIAM J. Sci. Comput., 30 (2008), pp. 290--311, https://doi.org/10.1137/060655742. 



\bibitem{HCEDJSAJW2014}
H. C. Elman, D. J. Silvester, and A. J. Wathen, Finite elements and fast iterative solvers: with applications in incompressible fluid dynamics, Oxford University Press, 2014.




\bibitem{CCPMAS1982}
C. C. Paige and M. A. Saunders, LSQR: An algorithm for sparse linear equations and sparse least squares, ACM Trans. Math. Software, 8 (1982), pp. 43--71, https://doi.org/10.1145/355984.355989. 




\bibitem{BFARDJSAJW1998}
B. Fischer, A. Ramage, D. J. Silvester, and A. J. Wathen, Minimum residual methods for augmented systems, BIT, 38 (1998), pp. 527--543, https://doi.org/10.1007/BF02510258.


\bibitem{SMA1995}
M. A. Saunders, Solution of sparse rectangular systems using LSQR and CRAIG, BIT, 35 (1995), pp. 588--604, https://doi.org/10.1007/BF01739829.


\bibitem{DCLFMAS2012}
D. C.-L. Fong, and M. A. Saunders, CG versus MINRES: An empirical comparison, SQU J. Science, 17 (2012), pp. 44--62, https://doi.org/10.24200/squjs.vol17iss1pp44-62.


\bibitem{DCLFMS2011}
D. C.-L. Fong, and M. Saunders, LSMR: an iterative algorithm for sparse least-quares problems, SIAM J. Sci. Comput., 33 (2011), pp. 2950--2971, https://doi.org/10.1137/10079687X.


\bibitem{CCPMAS1975}
C. C. Paige, and M. A. Saunders, Solution of sparse indefinite systems of linear equations, SIAM J. Numer. Anal., 12 (1975), pp. 617--629, https://doi.org/10.1137/0712047.


\bibitem{YSMHS1986}
Y. Saad, and M. H. Schultz, GMRES: A generalized minimal residual algorithm for solving nonsymmetric linear systems, SIAM J. Sci. Statist. Comput., 7 (1986), pp. 856--869, https://doi.org/10.1137/0907058.


\bibitem{HCEARDJS2014}
H. C. Elman, A. Ramage, and D. J. Silvester, IFISS: a computational laboratory for investigating incompressible flow problems, SIAM Rev., 56 (2014), pp. 261--273, https://doi.org/10.1137/120891393.


\bibitem{DSHEAR2014}
D. Silvester, H. Elman, and A. Ramage, Incompressible Flow and Iterative Solver Software (IFISS), Version 3.3 (2014), http://www.manchester.ac.uk/ifiss.


\bibitem{YS2003}
Y. Saad, Iterative Methods for Sparse Linear Systems, Second Ed., SIAM, 2003, https://doi.org/10.1137/1.9780898718003.


\bibitem{ECJLZS2024}
E. Carson, J. Liesen, and Z. Strako\v{s}, Towards understanding CG and GMRES through examples, Linear Algebra Appl., 692 (2024), pp. 241--291, https://doi.org/10.1016/j.laa.2024.04.003.

\bibitem{ADCKUR2024}
A. Dumitrasc, C. Kruse, and U. R\"{u}de, Deflation for the off-diagonal block in symmetric saddle point systems, SIAM J. Matrix Anal. Appl., 45 (2024), pp. 203--231, https://doi.org/10.1137/22M1537266.

%
%
%
%
%



\bibitem{MBGHGJL2005}
M. Benzi, G. H. Golub, and J. Liesen, Numerical solution of saddle point problems, Acta. Numer., 14 (2005), pp. 1--137, https://doi.org/10.1017/S0962492904000212.


\bibitem{MBMKN2006}
M. Benzi, and M. K. Ng, Preconditioned iterative methods for weighted Toeplitz least squares, SIAM J. Matrix Anal. Appl., 27 (2006), pp. 1106--1124, https://doi.org/10.1137/040616048.


\bibitem{MDAVDSDDS2010}
M. D'Apuzzo, V. De Simone, and D. di Serafino, On mutual impact of numerical linear algebra and large-scale optimization with focus on interior point methods, Comput. Optim. Appl., 45 (2010), pp. 283--310, https://doi.org/10.1007/s10589-008-9226-1. 


\bibitem{MPFDO2012}
M. P. Friedlander, and D. Orban, A primal-dual regularized interior-point method for convex quadratic problems, Math. Prog. Comp., 4 (2012), pp. 71--107, https://doi.org/10.1007/s12532-012-0035-2.


\bibitem{JPAJW2015}
J. Pestana, and A. J. Wathen, Natural preconditioning and iterative methods for saddle point systems, SIAM Rev., 57 (2015), pp. 51--71, https://doi.org/10.1137/130934921.

\bibitem{EVCMAO2000}
E. V. Chizhonkov, and M. A. Olshanskii, On the domain geometry dependence of the LBB condition, ESAIM Math. Model. Numer. Anal., 34 (2000), pp. 935--951, https://doi.org/10.1051/m2an:2000110.


%
%
%
%
%




\end{thebibliography}
\end{document}